\documentclass[a4paper,10pt]{amsart}

\usepackage{amssymb}
\usepackage[OT4]{fontenc}
\usepackage{texdraw}

\theoremstyle{definition}
\newtheorem{definition}{Definition}
\newtheorem{procedure}{Procedure}
\theoremstyle{plain}
\newtheorem{theorem}[definition]{Theorem}
\newtheorem{proposition}[definition]{Proposition}
\newtheorem{lemma}[definition]{Lemma}
\newtheorem{conjecture}[definition]{Conjecture}

\theoremstyle{remark}
\newtheorem{remark}[definition]{Remark}
\newtheorem{example}[definition]{Example}

\DeclareMathOperator{\edim}{edim}
\DeclareMathOperator{\vdim}{vdim}

\DeclareMathOperator{\mult}{mult}

\DeclareMathOperator{\red}{red}
\DeclareMathOperator{\symbred}{symbred}
\DeclareMathOperator{\symb}{symb}
\DeclareMathOperator{\Pic}{Pic}

\DeclareMathOperator{\rev}{rev}
\DeclareMathOperator{\supp}{supp}

\def\ass{\longleftarrow}
\def\field{\mathbb{K}}
\def\N{\mathbb{N}}

\def\Z{\mathbb{Z}}

\def\PP{\mathbb{P}}
\def\FF{\mathbb{F}}
\def\sys{\mathcal{L}}

\def\diag{\mathfrak{diag}}
\def\addarr{\stackrel{\text{add}}{\longrightarrow}}
\def\redarr{\stackrel{\text{red}}{\longrightarrow}}
\def\symbredarr{\stackrel{\text{symbred}}{\longrightarrow}}
\def\symbx{\mathbf{x}}
\def\ldot{\cdot}
\def\up{\uparrow}
\def\lineq{\sim}

\let\to\longrightarrow
\let\mapsto\longmapsto

\begin{document}

\title{Special homogeneous linear systems on Hirzebruch surfaces}

\author{Marcin Dumnicki}

\dedicatory{
Institute of Mathematics, Jagiellonian University, \\
ul. \L{}ojasiewicza 6, 30-348 Krak\'ow, Poland \\
Email address: Marcin.Dumnicki@im.uj.edu.pl\\
}



\begin{abstract}
The Segre-Gimigliano-Harbourne-Hirschowitz Conjecture can be
naturally formulated for Hirzebruch surfaces $\FF_n$. We show
that this Conjecture holds for imposed base points of equal
multiplicity bounded by 8.
\end{abstract}

\maketitle

\section{Linear systems on Hirzebruch surfaces}

Our goal is to prove Conjecture \ref{mainconj} for 
linear systems on Hirzebruch surfaces with imposed base points
of equal multiplicity bounded by $8$. This Conjecture,
being a natural reformulation of the Segre-Harbourne-Gimigliano-Hirschowitz
Conjecture,
has been stated in \cite[Conjecture 2.6]{Laf}. In the same paper it is shown
(Theorem 7.1) that this Conjecture holds for systems with
imposed base points of equal multiplicity bounded by $3$.
We will also give another proof of \cite[Proposition 2.7]{Laf}, where
the proof contains a serious mistake (for more details see the proof
of Proposition \ref{easyfn}).

Our method will also work for greater values of multiplicities, but
the computational part (realized with the help of computers) becomes very
large and time-consuming. But it is possible to carry our computations
further to obtain the proof for $m_1=\dots=m_r=9,10,\dots$ or to find a
counterexample.

The author would like to thank Micha\l{} Kapustka and
Tomasz Szemberg for valuable discussions.

By $\FF_n$, $n \geq 0$, we denote the rational ruled surface (called the \emph{$n$-th
Hirzebruch surface}) given by
$\FF_n = \PP(\mathcal O_{\PP^1} \oplus \mathcal O_{\PP^1}(n))$
over the field $\field$ of characteristic $0$.
The Picard group $\Pic(\FF_n)$ can be freely generated by the class of a fiber $F_n$
and the class of the section $H_n$ such that $F_n^2=0$, $H_n^2=n$, $F_n\ldot H_n=1$.
The irreducible section with self-intersection $-n$ will be denoted
by $\Gamma_n$, we have $\Gamma_n \in |H_n-nF_n|$. The class
of $\Gamma_n$ in $\Pic(\FF_n)$ will also be denoted by $\Gamma_n$. Let $a$, $b$ be
integers. By $\sys_n(a,b)$ we will denote the complete linear system
associated to the line bundle $aF_n+bH_n$.

\begin{lemma}
If on $\FF_n$ the class $aF_n+bH_n$ contains an effective divisor 
then there exists non-negative integers
$a'$, $b'$, $q$ ($q > 0$ if and only if $a<0$) such that 
the base locus of $|aF_n+bH_n|$ is $q\Gamma_n$
and $aF_n+bH_n$ is linearly equivalent to $q\Gamma_n+a'F_n+b'H_n$. Moreover, we have
$$\dim \sys_n(a,b) = \frac{(b'+1)(2a'+2+nb')}{2}-1$$
\end{lemma}

\begin{proof}
For the proof see \cite[Proposition 2.2]{Laf}.
\end{proof}

Now we pick $r$ points $p_1,\dots,p_r \in \FF_n$ in general position,
let $m_1,\dots,m_r$ be non-negative integers.
By $\sys_n(a,b;m_1,\dots,m_r)$ we denote the linear system
of curves in $\sys_n(a,b)$ passing through points $p_1,\dots,p_r$
with multiplicities at least $m_1,\dots,m_r$, respectively.
The points $p_1,\dots,p_r$ will be called \emph{imposed base points}.
The dimension of this system will be denoted by
$\dim \sys_n(a,b;m_1,\dots,m_r)$. Define the
\emph{virtual dimension}
\begin{align*}
\vdim \sys_n(a,b;m_1,\dots,m_r) & = \dim \sys_n(a,b) - \sum_{j=1}^{r} \binom{m_j+1}{2} \\
\intertext{and the \emph{expected dimension}}
\edim \sys_n(a,b;m_1,\dots,m_r) & = \max \{ \vdim \sys_n(a,b;m_1,\dots,m_r), -1 \}.\\
\intertext{We have}
\dim \sys_n(a,b;m_1,\dots,m_r) & \geq \edim \sys_n(a,b;m_1,\dots,m_r).
\end{align*}
If this inequality is strict then the system $\sys_n(a,b;m_1,\dots,m_r)$
is said to be \emph{special}, \emph{non-special} otherwise.
The system of negative dimension will be called \emph{empty}.

A natural question is: when a given system is special, and if there exists
a geometric explanation to the non-speciality. This can be done by
considering $-1$-systems.

To introduce the notion of $-1$-system and $-1$-speciality define
the \emph{intersection number} of
$L=\sys_n(a,b;m_1,\dots,m_r)$ and $L'=\sys_n(a,b;m_1',\dots,m_r')$
\begin{align*}
L \ldot L' & = (aF_n+bH_n) \ldot (a'F_n+b'H_n) - \sum_{j=1}^{r} m_jm_j' \\
& = ab'+a'b+nbb' - \sum_{j=1}^{r} m_jm_j'.
\end{align*}
Observe that if, for nonempty systems,
$L \ldot L' < 0$ then
these systems must have a common component. The intersection number of
two systems $L$ and $L'$ can also be defined by taking the blow-up $\pi:S \to \FF_n$ at
imposed base points and putting $L \ldot L' = \pi^{\star}(L) \ldot \pi^{\star}(L')$
on $S$ (see \cite{Laf}).

\begin{definition}
The system $E=\sys_n(a,b;m_1,\dots,m_r)$, $a,b \geq 0$, satisfying
$$\dim E = \vdim E = 0, \quad E \ldot E = -1,$$
with irreducible member is called \emph{$-1$-system}.
\end{definition}

\begin{procedure}
\label{proccr}
Let $L=\sys_n(a,b;m_1,\dots,m_r)$, consider the following procedure:
\begin{itemize}
\item Step 1. Put $M \ass L$.
\item Step 2. If $M \ldot \Gamma_n < 0$ then take $M \ass M-\Gamma_n$
and go back to Step 2.
\item Step 3. If $M \ldot E < 0$ for some $-1$-system $E$ then
take $M \ass M-E$ and go back to Step 2.
\end{itemize}
The procedure terminates after a finite number of steps.
\end{procedure}

\begin{definition}
If, for $M$ and $L$ as above, $\edim M > \edim L$ then $L$ will be called \emph{$-1$-special}.
\end{definition}

Observe that if $L$ is $-1$-special then
$$\dim L = \dim M \geq \edim M > \edim L,$$
so $L$ is special. In \cite[Conjecture 2.6]{Laf} it is conjectured that

\begin{conjecture}
\label{mainconj}
The system $\sys_n(a,b;m_1,\dots,m_r)$ is special if and only if it is $-1$-special.
\end{conjecture}

The analogous Conjecture for the projective plane was stated by several
authors and is known as Segre-Harbourne-Gimigliano-Hirschowitz Conjecture (SHGH for short).
More on this Conjecture can be found e.g. in \cite{homo10}, some recent
results are listed in \cite{Dumquasi}.

Since we are interested mainly in homogeneous systems,
we will use the notation $m^{\times r}$ for repeated multiplicities.

\begin{example}
Let us consider $\sys_6(0,4;3^{\times 11})$. Observe that
$\Gamma_n \in \sys_n(-n,1)$. We have 
$$\sys_6(0,4;3^{\times 11}) \ldot \sys_6(-6,1) = 0,$$
so we pass to Step 3 in Procedure \ref{proccr}. For the $-1$-system $E=\sys_6(2,1;1^{\times 11})$
we have 
$$\sys_6(0,4;3^{\times 11}) \ldot E = 8+24-33=-1,$$
so we must take new system $\sys_6(-2,3;2^{\times 11})$. In Procedure \ref{proccr}, Step 2
$$\sys_6(-2,3;2^{\times 11}) \ldot \Gamma_6 = -2,$$
hence we take out the $-n$-section from the base locus and
obtain $\sys_6(4,2;2^{\times 11})$, which is equal to $2E$.
Consequently we have that $\sys_6(0,4;3^{\times 11}) = \Gamma_6+3E$,
which is non-empty, and since
$\vdim \sys_6(0,4;3^{\times 11}) = -2$,
it is $-1$-special.
\end{example}

\section{Linear systems over $\PP^2$}

\begin{definition}
Let $d,m_1,\dots,m_r,k_1,\dots,k_s$ be non-negative integers.
Pick a general line $\ell \subset \PP^2$, pick
points $p_1,\dots,p_r$ in general position, pick points
$q_1,\dots,q_s \in \ell$ also in general position on the line $\ell$.
By
$$\sys(d;m_1,\dots,m_r,\overline{k_1,\dots,k_s})$$
we denote the linear system of curves in $\PP^2$ of degree $d$ with
multiplicities at least $m_1,\dots,m_r$, $k_1,\dots,k_s$ at
$p_1,\dots,p_r$, $q_1,\dots,q_s$ respectively.
The dimension of this system will be denoted by
$$\dim \sys(d;m_1,\dots,m_r,\overline{k_1,\dots,k_s}).$$
Define the
\emph{virtual dimension}
\begin{align*}
\vdim \sys(d;m_1,\dots,m_r,\overline{k_1,\dots,k_s}) & = \binom{d+2}{2} - 1 - \sum_{j=1}^{r} \binom{m_j+1}{2} - \sum_{j=1}^{s} \binom{k_j+1}{2} \\
\intertext{and the \emph{expected dimension}}
\edim \sys(d;m_1,\dots,m_r,\overline{k_1,\dots,k_s}) & = \max \{ \vdim \sys(d;m_1,\dots,m_r,\overline{k_1,\dots,k_s}), -1 \}.\\
\intertext{We have}
\dim \sys(d;m_1,\dots,m_r,\overline{k_1,\dots,k_s}) & \geq \edim \sys(d;m_1,\dots,m_r,\overline{k_1,\dots,k_s}).
\end{align*}
If this inequality is strict then the system $\sys(d;m_1,\dots,m_r,\overline{k_1,\dots,k_s})$
is said to be \emph{special}, \emph{non-special} otherwise.
We also have the \emph{intersection number}
\begin{multline*}
\sys(d;m_1,\dots,m_r,\overline{k_1,\dots,k_s}) \ldot
\sys(d';m_1',\dots,m_r',\overline{k_1',\dots,k_s'}) \\
= dd' - \sum_{j=1}^{r} m_jm_j' - \sum_{j=1}^{s} k_jk_j'.
\end{multline*}
\end{definition}

Again, we define $-1$-system and $-1$-speciality.

\begin{definition}
The system $E=\sys(d;m_1,\dots,m_r,\overline{k_1,\dots,k_s})$ satisfying
$$\dim E = \vdim E = 0, \quad E \ldot E = -1,$$
with irreducible member is called \emph{$-1$-system}.
\end{definition}

\begin{procedure}
\label{proccrp2}
Let $L=\sys(d;m_1,\dots,m_r,\overline{k_1,\dots,k_s})$, consider the following procedure:
\begin{itemize}
\item Step 1. Put $M \ass L$.
\item Step 2. If $M \ldot \sys(1;\overline{1^{\times s}}) < 0$ then take $M \ass M-\sys(1;\overline{1^{\times s}})$
and go back to Step 2.
\item Step 3. If $M \ldot E < 0$ for some planar $-1$-system $E$ then
take $M \ass M-E$ and go back to Step 2.
\end{itemize}

Step 2 should be understood as follows: for
$M=\sys(d';m_1',\dots,m_r',\overline{k_1',\dots,k_s'})$ we have
$$M \ldot \sys(1;\overline{1^{\times s}}) = d-\sum_{j=1}^{s} k_j$$
and if this number is negative then the line $\ell$ lies in the base locus
of $M$ and can be taken out as follows:
$$M - \sys(1;\overline{1^{\times s}}) = \sys(d-1;m_1,\dots,m_r,\overline{k_1-1,\dots,k_s-1}).$$
\end{procedure}

\begin{definition}
If, after Procedure \ref{proccrp2} terminates, $\edim M > \edim L$ then $L$ will be called \emph{$-1$-special}.
\end{definition}

\begin{example}
\label{ex28}
Let us consider $L=\sys(28;24,3^{\times 11},\overline{4^{\times 7}})$.
In the next section we will show that $\dim L = \dim \sys_6(0,4;3^{\times 11})$
and $L$ is $-1$-special (resp. special) if and only if $\sys_6(0,4;3^{\times 11})$
is $-1$-special (resp. special).
We have $L \ldot E = -1$ for $E = \sys(9;8,1^{\times 11},\overline{1^{\times 7}})$.
The residual system $L-E = \sys(19;16,2^{\times 11},\overline{3^{\times 7}})$
has the line in the base locus. Continuing this way we will have
$L = \sys(1;\overline{1^{\times 7}}) + 3E$, and since $E$ is a $-1$-system,
$L$ is $-1$-special.
\end{example}

\begin{remark}
The original SHGH Conjecture states
that for a plane system with imposed base points in general position
the speciality is equivalent
to the $-1$-speciality. For a system with collinear imposed base points 
it is natural to extend the definition of the $-1$-speciality as above,
which should be called the \emph{negative speciality}, since the self-intersection
of the line passing through $s$ imposed base points is equal to $1-s$.
\end{remark}

We will often consider plane systems and systems on Hirzebruch surfaces.
Therefore, we will consequently use the notation: $\sys$ without number
denotes always the system on $\PP^2$, while $\sys_n$ the system on $\FF_n$.

\begin{proposition}
\label{cremona}
Let $d,m_1,\dots,m_r,k_1\dots,k_s$ be non-negative integers.
\begin{itemize}
\item
Let $k=d-m_1-m_2-m_3$, let
$m_j^{\star}=\max\{m_j+k,0\}$ for $j=1,2,3$.
Then
$$\dim \sys(d;m_1,\dots,m_r) = \dim \sys(d+k;m_1^{\star},m_2^{\star},m_3^{\star},m_4,\dots,m_r).$$
\item
Let $k=d-m_1-m_2-k_1$, let $m_1^{\star}=\max\{m_1+k,0\}$, $m_2^{\star}=\max\{m_2+k,0\}$,
$k_1^{\star}=\max\{k_1+k,0\}$.
Then
\begin{multline*}
\dim \sys(d;m_1,\dots,m_r,\overline{k_1,\dots,k_s}) \\
= \dim \sys(d+k;m_1^{\star},m_2^{\star},m_3,\dots,m_r,\overline{k_1^{\star},k_2,\dots,k_s}).
\end{multline*}
\end{itemize}
Moreover, let $L$ denote the original system and $L^{\star}$ 
the system after transformation. Then either $\edim L = \edim L^{\star}$
or $\edim L^{\star} > \edim L$ and $L$ is $-1$-special.
\end{proposition}

\begin{proof}
Let $L=\sys(d;m_1,\dots,m_r)$, let $L^{\star}=\sys(d+k;m_1^{\star},m_2^{\star},m_3^{\star},m_4,\dots,m_r)$.
To show the first equality we must check if $\sys(1;1,1,0)$, $\sys(1;1,0,1)$
or $\sys(1;0,1,1)$ are in the base locus of $L$ and write
$$L = q_1\sys(1;1,1,0)+q_2\sys(1;1,0,1)+q_3\sys(1;0,1,1)+\widetilde{L}$$
for the system $\widetilde{L}$ without these lines in the base locus. It follows
that $\dim L = \dim \widetilde{L}$ and if $\edim \widetilde{L} > \edim L$
then $L$ is $-1$-special. To complete the proof
observe that applying the standard birational transformation
(so called \emph{Cremona transformation}) based on points $p_1$, $p_2$, $p_3$
to $\widetilde{L}$ we obtain the system $L^{\star}$. By a simple calculation
we can show that $\edim \widetilde{L} = \edim L^{\star}$.

To see the second
equality observe that $\sys(1;1,0,0)$ is invariant under Cremona transformation,
so the line passing through exactly one of the three points will be preserved.
\end{proof}

\begin{remark}
We can apply the above to any three multiplicities, since we can
permute imposed points.
\end{remark}

\begin{example}
Let us again (see Example \ref{ex28}) consider $L=\sys(28;24,3^{\times 11},\overline{4^{\times 7}})$.
This time we are only interested in showing that $L$ is non-empty.
We can make Cremona transformation based on points with multiplicity
$24$, $3$ and $4$ to obtain
$\sys(25;21,3^{\times 10},\overline{4^{\times 6},1}).$
We can repeat this $6$ more times, which leads us to a system
$\sys(7;3,3^{\times 4},\overline{1^{\times 7}})$. If
$\sys(6;3^{\times 5})$ is non-empty then $L$ will also be non-empty.
On one hand $\vdim \sys(6;3^{\times 5}) = -3$, but on the other hand, applying Cremona,
we have that $\dim \sys(6;3^{\times 5}) = \dim \sys(3;3,3) = 0$.
\end{example}

\section{From Hirzebruch surface to $\PP^2$}

\begin{proposition}
\label{fromhitop2}
Let $n \geq 0$. For any non-negative integers $a$, $b$, $a'$, $b'$,
$m_1,\dots,m_r$, $m_1',\dots,m_r'$ we have
\begin{gather*}
\dim \sys_n(a,b;m_1,\dots,m_r) = \dim \sys(a+(n+1)b;a+nb,m_1,\dots,m_r,\overline{b^{\times (n+1)}}), \\
\vdim \sys_n(a,b;m_1,\dots,m_r) = \vdim \sys(a+(n+1)b;a+nb,m_1,\dots,m_r,\overline{b^{\times (n+1)}}), \\
\sys_n(a,b;m_1,\dots,m_r) \ldot \sys_n(a',b';m_1',\dots,m_r') \qquad \qquad \qquad \qquad \qquad \qquad \qquad \\ 
= \sys(a+(n+1)b;a+nb,m_1,\dots,m_r,\overline{b^{\times (n+1)}}) \\
\qquad \qquad \qquad \qquad \ldot \sys(a'+(n+1)b';a'+nb',m_1',\dots,m_r',\overline{b'^{\times (n+1)}}).
\end{gather*}
\end{proposition}

\begin{proof}
We will use $\lineq$ for linear equivalence of divisors.
By a straightforward calculation we show the above for
the virtual dimension and the intersection number.

Let $L=\sys_n(a,b;m_1,\dots,m_r)$. For $n \geq 2$ consider the blow-up $\pi_1:X \to \FF_n$ 
of a point $p$ on a fiber $F_p$, $p \notin \Gamma_n$. We can also assume
that there are no imposed base point on $F_p$.
Let $E$ be the exceptional divisor of $\pi_1$,
let $\pi_1^{\star}(F_p)=\widetilde{F_p}+E$, $\pi_1^{\star}(\Gamma_n)=\widetilde{\Gamma}$.
Now blow down $\widetilde{F_p}$ (which has self-intersection equal to $-1$)
with $\pi_2:X \to \FF_{n-1}$ (see Figure \ref{blow1}), let $q=\pi_2(\widetilde{F_p})$,
let $F_q$ be the fiber on $\FF_{n-1}$ passing through $q$.
The above is often called an {\it elementary transformation}.
\begin{figure}[ht!]
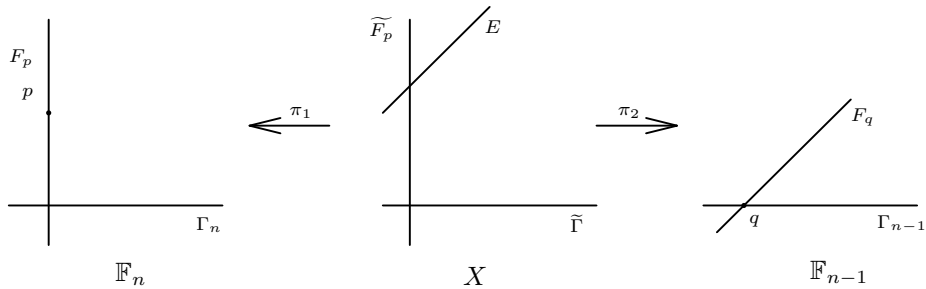

\begin{texdraw}
\drawdim pt
\arrowheadtype t:V
\move(0 0)
\lvec(80 0)
\move(15 -15)
\lvec(15 70)
\move(15 35)
\fcir f:0 r:1
\htext(70 -10){{\scriptsize $\Gamma_n$}}
\htext(0 52){{\scriptsize $F_p$}}
\htext(5 40){{\scriptsize $p$}}

\move(120 30)
\avec(90 30)
\htext(105 33){{\scriptsize $\pi_1$}}

\move(140 0)
\lvec(220 0)
\move(150 -15)
\lvec(150 70)
\move(140 35)
\lvec(180 75)
\htext(210 -10){{\scriptsize $\widetilde{\Gamma}$}}
\htext(135 62){{\scriptsize $\widetilde{F_p}$}}
\htext(178 65){{\scriptsize $E$}}

\move(220 30)
\avec(250 30)
\htext(228 33){{\scriptsize $\pi_2$}}

\move(260 0)
\lvec(340 0)
\move(265 -10)
\lvec(315 40)
\move(275 0)
\fcir f:0 r:1
\htext(325 -10){{\scriptsize $\Gamma_{n-1}$}}
\htext(315 30){{\scriptsize $F_q$}}
\htext(277 -8){{\scriptsize $q$}}

\htext(40 -30){$\FF_n$}
\htext(170 -30){$X$}
\htext(300 -30){$\FF_{n-1}$}
\end{texdraw}
\caption{Elementary transformation between $\FF_n$ and $\FF_{n-1}$}\label{blow1}
\end{figure}
Let $F$ denote the class of $\widetilde{F_p}$ in $\Pic(X)$,
we will denote the classes of $\widetilde{\Gamma}$ and $E$ by $\Gamma$ and $E$,
respectively. We have
$\pi_2^{\star}(F_n) \lineq F+E$, $\pi_2^{\star}(\Gamma_{n-1})=\Gamma+F$. Moreover,
$\pi_1^{\star}(H_n) \lineq \pi_1^{\star}(\Gamma_n+nF_n) \lineq \Gamma+nF+nE$,
$\pi_2^{\star}(H_{n-1}) \lineq \pi_2^{\star}(\Gamma_{n-1}+(n-1)F_{n-1})=\Gamma+nF+(n-1)E$.
The strict transform of the class of a curve in $\sys_n(a,b)$ is 
\begin{align*}
a\pi_1^{\star}(F_n)+b\pi_1^{\star}(H_n) & \lineq a(F+E)+b\Gamma+nbF+nbE \\
& \lineq b\Gamma+nbF+b(n-1)E+(a+b)(F+E)-bF \\
& \lineq (a+b)\pi_2^{\star}(F_{n-1})+b\pi_2^{\star}(H_{n-1})-bF.
\end{align*}
After blowing a curve from $\sys_n(a,b)$ up by $\pi_1$ and blowing down by $\pi_2$
we obtain the curve belonging to
$\sys_{n-1}(a+b,b;b)$, and the point with multiplicity $b$
lies on $\Gamma_{n-1}$. Since $E$
does not belong to the base locus of $\pi_1^{\star}(aF_n+bH_n)$, and $\widetilde{F_p}$ belongs
to the base locus exactly $b$ times, from Leray spectral sequence
we have
$$\dim \sys_n(a,b;m_1,\dots,m_r) = \dim \sys_{n-1}(a+b,b;b,m_1,\dots,m_r).$$

By repeating the above process we will end up with the system
$$\sys_1(a+(n-1)b,b;b^{\times (n-1)},m_1,\dots,m_r),$$
where $n-1$ imposed
base points lies generically on the $-1$-curve $\Gamma_1$.

The surface $\FF_1$ is isomorphic to $\PP^2$ blown up in one point with the exceptional divisor
$\Gamma_1$. Take fibers $F_p$ and $F_{p'}$ passing through general points
$p$ and $p'$, respectively. Let $H_{pp'} \subset \FF_1$ be the 
strict transform of the line (from $\PP^2$) joining $p$ and $p'$.
Let $\pi_1: X \to \FF_1$ be the sequence of two blow-ups:
of $p$ and $p'$ with exceptional divisors $E$ and $E'$, respectively.
Let $\pi_1^{\star}(F_p)=\widetilde{F_p}+E$, $\pi_1^{\star}(F_{p'})=\widetilde{F_{p'}}+E'$,
$\pi_1^{\star}(\Gamma_1)=\widetilde{\Gamma}$, $\pi_1^{\star}(H_{pp'})=\widetilde{H_{pp'}}+E+E'$.
Let $F$ (resp. $F'$, $\Gamma$) denote the class of $\widetilde{F_p}$ (resp.
$\widetilde{F_{p'}}$, $\widetilde{\Gamma}$) in $\Pic(X)$.
We have $\pi_1^{\star}(H_1) \lineq \Gamma+E+F \lineq \Gamma+E'+F'$.
Now take the sequence of three blow-downs $\pi_2:X \to Y$:
of $\widetilde{F_p}$, $\widetilde{F_{p'}}$ and $\widetilde{H_{pp'}}$ (see Figure \ref{blow2}).
\begin{figure}[ht!]
\begin{texdraw}
\drawdim pt
\arrowheadtype t:V
\move(0 0)
\lvec(80 0)
\move(15 -15)
\lvec(15 60)
\move(60 -15)
\lvec(60 60)
\move(0 35)
\lvec(80 35)
\move(15 35)
\fcir f:0 r:1
\move(60 35)
\fcir f:0 r:1
\htext(70 -10){{\scriptsize $\Gamma_1$}}
\htext(2 52){{\scriptsize $F_p$}}
\htext(45 52){{\scriptsize $F_{p'}$}}
\htext(7 37){{\scriptsize $p$}}
\htext(51 37){{\scriptsize $p'$}}
\htext(70 38){{\scriptsize $H_{pp'}$}}

\move(120 30)
\avec(90 30)
\htext(105 33){{\scriptsize $\pi_1$}}

\move(140 0)
\lvec(220 0)
\move(150 -15)
\lvec(150 60)
\move(195 -15)
\lvec(195 60)
\move(140 35)
\lvec(190 85)
\move(185 35)
\lvec(235 85)
\move(150 75)
\lvec(240 75)
\htext(210 -10){{\scriptsize $\widetilde{\Gamma}$}}
\htext(140 52){{\scriptsize $\widetilde{F_p}$}}
\htext(182 52){{\scriptsize $\widetilde{F_{p'}}$}}
\htext(178 65){{\scriptsize $E$}}
\htext(223 65){{\scriptsize $E'$}}
\htext(150 78){{\scriptsize $\widetilde{H_{pp'}}$}}

\move(220 30)
\avec(250 30)
\htext(228 33){{\scriptsize $\pi_2$}}

\move(260 0)
\lvec(340 0)
\move(265 -10)
\lvec(310 35)
\move(285 35)
\lvec(330 -10)
\move(275 0)
\fcir f:0 r:1
\move(320 0)
\fcir f:0 r:1
\move(297.5 22.5)
\fcir f:0 r:1
\htext(335 -10){{\scriptsize $\overline{\Gamma}$}}
\htext(275 12){{\scriptsize $\overline{E}$}}
\htext(310 12){{\scriptsize $\overline{E'}$}}

\htext(32 -30){$\FF_1$}
\htext(163 -30){$X$}
\htext(295 -30){$\PP^2$}
\end{texdraw}
\caption{Transformation between $\FF_1$ and $\PP^2$}\label{blow2}
\end{figure}
The above rational transformation $\FF_1 \to Y$ is nothing else than realizing the Cremona transformation
of $\PP^2$ 
by three blow-ups and three blow-downs, but we start with one point
blown-up already. It follows that $Y=\PP^2$ and
$\pi_2^{\star}(\overline{\Gamma}) \lineq \Gamma+F+F'$ for
the class of
the line $\overline{\Gamma} = \pi_2(\widetilde{\Gamma})$ in $\PP^2$.
Now, for $\sys_1(a,b)$, we have
\begin{align*}
a\pi_1^{\star}(F_1)+b\pi_1^{\star}(H_1) & \lineq a(F+E)+b(\Gamma+E+F) \\
& \lineq (a+2b)(\Gamma+F+F')-(a+b)(\Gamma+F'-E)-bF-bF' \\
& \lineq (a+2b)\pi_2^{\star}(\overline{\Gamma})-(a+b)H-bF-bF'.
\end{align*}
After blowing down we obtain the curve in $\sys(a+2b;a+b,b,b)$.
The section $\Gamma_1$ is preserved and mapped to the line
$\overline{\Gamma}$. In consequence we have that, taking $p_1,\dots,p_r \in \FF_n$
in general position,
$$
\dim \sys_n(a,b;m_1,\dots,m_r) = \dim \sys(a+(n+1)b;a+nb,m_1,\dots,m_r,\overline{b^{\times n+1}}).
$$
For $\FF_0=\PP^1\times \PP^1$ the easy proof is left to the reader.
\end{proof}

\begin{remark}
Observe that $-1$-systems on $\FF_n$ are transformed into $-1$-systems on $\PP^2$,
since the dimension, virtual dimension, self-intersection and irreducibility
is preserved. The section $\Gamma_n$ is mapped into a line which contains
$n+1$ imposed base points with multiplicities $b^{\times (n+1)}$. This means,
in particular, that the system $\sys_n(a,b;m_1,\dots,m_r)$
is $-1$-special if and only if the planar system 
$$\sys(a+(n+1)b;a+nb,m_1,\dots,m_r,\overline{b^{\times (n+1)}})$$
is $-1$-special (compare Procedures \ref{proccr} and \ref{proccrp2}). Additionally we will see
that the dimension of each considered $-1$-system 
$$\sys(d;m_1,\dots,m_r,\overline{k_1,\dots,k_s})$$
on $\PP^2$
will remain $0$ after assigning all base points in general position,
i.e. 
$$\sys(d;m_1,\dots,m_r,k_1,\dots,k_s)$$
will also be a $-1$-system
(see the last section). Therefore we can state the following Conjecture.
\end{remark}

\begin{conjecture}
\begin{multline*}
\dim \sys(d;m_1,\dots,m_r,\overline{k_1,\dots,k_s}) \\
= \max_{j=0,\dots,d} \dim \sys(d-j;m_1,\dots,m_r,\max \{k_1-j,0\},\dots,\max \{k_s-j,0\}).
\end{multline*}
\end{conjecture}

\section{Diagrams and reductions}

\begin{definition}
Let $a_1,\dots,a_s$ be non-negative integers. Set
$a_{j}=0$ for $j > s$ and define the \emph{diagram}
$$\diag(a_1,\dots,a_s) = \{ (x,y) \in \N^2 : y < a_{x+1} \}.$$
We will also write $[a]^{\times p}$ for
$$\underbrace{a,\dots,a}_{p},$$
and $\diag(a_1,\dots,a_s)+\diag(b_1,\dots,b_p)$ for $\diag(a_1,\dots,a_s,b_1,\dots,b_p)$.
\end{definition}

\begin{example}
$$
\begin{array}{ccc}

\begin{texdraw}
\drawdim pt
\move(0 0)
\fcir f:0 r:0.5
\move(0 10)
\fcir f:0 r:0.5
\move(0 20)
\fcir f:0 r:0.5
\move(0 30)
\fcir f:0 r:0.5
\move(0 40)
\fcir f:0 r:0.5
\move(10 0)
\fcir f:0 r:0.5
\move(10 10)
\fcir f:0 r:0.5
\move(10 20)
\fcir f:0 r:0.5
\move(10 30)
\fcir f:0 r:0.5
\move(10 40)
\fcir f:0 r:0.5
\move(20 0)
\fcir f:0 r:0.5
\move(20 10)
\fcir f:0 r:0.5
\move(20 20)
\fcir f:0 r:0.5
\move(20 30)
\fcir f:0 r:0.5
\move(20 40)
\fcir f:0 r:0.5
\move(30 0)
\fcir f:0 r:0.5
\move(30 10)
\fcir f:0 r:0.5
\move(30 20)
\fcir f:0 r:0.5
\move(30 30)
\fcir f:0 r:0.5
\move(30 40)
\fcir f:0 r:0.5
\move(40 0)
\fcir f:0 r:0.5
\move(40 10)
\fcir f:0 r:0.5
\move(40 20)
\fcir f:0 r:0.5
\move(40 30)
\fcir f:0 r:0.5
\move(40 40)
\fcir f:0 r:0.5
\move(50 0)
\fcir f:0 r:0.5
\move(50 10)
\fcir f:0 r:0.5
\move(50 20)
\fcir f:0 r:0.5
\move(50 30)
\fcir f:0 r:0.5
\move(50 40)
\fcir f:0 r:0.5
\move(60 0)
\fcir f:0 r:0.5
\move(60 10)
\fcir f:0 r:0.5
\move(60 20)
\fcir f:0 r:0.5
\move(60 30)
\fcir f:0 r:0.5
\move(60 40)
\fcir f:0 r:0.5
\arrowheadtype t:V
\move(0 0)
\avec(90 0)
\move(0 0)
\avec(0 70)
\htext(96 0){$\mathbb{N}$}
\htext(-13 50){$\mathbb{N}$}
\move(0 0)
\fcir f:0 r:1.5
\move(10 0)
\fcir f:0 r:1.5
\move(10 10)
\fcir f:0 r:1.5
\move(20 0)
\fcir f:0 r:1.5
\move(20 10)
\fcir f:0 r:1.5
\move(20 20)
\fcir f:0 r:1.5
\move(30 0)
\fcir f:0 r:1.5
\move(30 10)
\fcir f:0 r:1.5
\move(30 20)
\fcir f:0 r:1.5
\move(30 30)
\fcir f:0 r:1.5
\move(40 0)
\fcir f:0 r:1.5
\move(40 10)
\fcir f:0 r:1.5
\move(40 20)
\fcir f:0 r:1.5
\move(50 0)
\fcir f:0 r:1.5
\move(50 10)
\fcir f:0 r:1.5
\textref h:C v:C
\htext(40 -15){diagram $\diag(1,2,3,4,3,2)$}
\end{texdraw}

&

\hspace{1cm}

&

\begin{texdraw}
\drawdim pt
\move(0 0)
\fcir f:0 r:0.5
\move(0 10)
\fcir f:0 r:0.5
\move(0 20)
\fcir f:0 r:0.5
\move(0 30)
\fcir f:0 r:0.5
\move(0 40)
\fcir f:0 r:0.5
\move(0 50)
\fcir f:0 r:0.5
\move(10 0)
\fcir f:0 r:0.5
\move(10 10)
\fcir f:0 r:0.5
\move(10 20)
\fcir f:0 r:0.5
\move(10 30)
\fcir f:0 r:0.5
\move(10 40)
\fcir f:0 r:0.5
\move(10 50)
\fcir f:0 r:0.5
\move(20 0)
\fcir f:0 r:0.5
\move(20 10)
\fcir f:0 r:0.5
\move(20 20)
\fcir f:0 r:0.5
\move(20 30)
\fcir f:0 r:0.5
\move(20 40)
\fcir f:0 r:0.5
\move(20 50)
\fcir f:0 r:0.5
\move(30 0)
\fcir f:0 r:0.5
\move(30 10)
\fcir f:0 r:0.5
\move(30 20)
\fcir f:0 r:0.5
\move(30 30)
\fcir f:0 r:0.5
\move(30 40)
\fcir f:0 r:0.5
\move(30 50)
\fcir f:0 r:0.5
\move(40 0)
\fcir f:0 r:0.5
\move(40 10)
\fcir f:0 r:0.5
\move(40 20)
\fcir f:0 r:0.5
\move(40 30)
\fcir f:0 r:0.5
\move(40 40)
\fcir f:0 r:0.5
\move(40 50)
\fcir f:0 r:0.5
\move(50 0)
\fcir f:0 r:0.5
\move(50 10)
\fcir f:0 r:0.5
\move(50 20)
\fcir f:0 r:0.5
\move(50 30)
\fcir f:0 r:0.5
\move(50 40)
\fcir f:0 r:0.5
\move(50 50)
\fcir f:0 r:0.5
\arrowheadtype t:V
\move(0 0)
\avec(80 0)
\move(0 0)
\avec(0 80)
\htext(86 0){$\mathbb{N}$}
\htext(-13 60){$\mathbb{N}$}
\move(0 0)
\fcir f:0 r:1.5
\move(0 10)
\fcir f:0 r:1.5
\move(0 20)
\fcir f:0 r:1.5
\move(0 30)
\fcir f:0 r:1.5
\move(0 40)
\fcir f:0 r:1.5
\move(10 0)
\fcir f:0 r:1.5
\move(10 10)
\fcir f:0 r:1.5
\move(10 20)
\fcir f:0 r:1.5
\move(10 30)
\fcir f:0 r:1.5
\move(10 40)
\fcir f:0 r:1.5
\move(20 0)
\fcir f:0 r:1.5
\move(20 10)
\fcir f:0 r:1.5
\move(20 20)
\fcir f:0 r:1.5
\move(20 30)
\fcir f:0 r:1.5
\move(20 40)
\fcir f:0 r:1.5
\move(30 0)
\fcir f:0 r:1.5
\move(30 10)
\fcir f:0 r:1.5
\move(30 20)
\fcir f:0 r:1.5
\move(40 0)
\fcir f:0 r:1.5
\move(40 10)
\fcir f:0 r:1.5
\move(40 20)
\fcir f:0 r:1.5
\textref h:C v:C
\htext(35 -15){diagram $\diag([5]^{\times 3},[3]^{\times 2})$}
\end{texdraw}

\end{array}
$$
\end{example}

\begin{definition}
Let $m > 0$, let $a_1,\dots,a_s$, $b_1,\dots,b_m$ be non-negative integers,
$b_j > 0$ for $j=1,\dots,m$.
Define $r_1,\dots,r_m$ inductively (beginning with $r_m$) to be
$$
r_j = 
\begin{cases}
b_j & \text{if } b_j < m, \\
\max ( \{1,\dots,m\} \setminus \{r_{j+1},\dots,r_m\} ) & \text{if } b_j \geq m.
\end{cases}
$$
If $\{r_1,\dots,r_m\}=\{1,\dots,m\}$ then we say that $D=\diag(a_1,\dots,a_s,b_1,\dots,b_m)$
is \emph{$m$-reducible} and define the \emph{$m$-reduction of $D$}
$$
\red_m(D) = \diag(a_1,\dots,a_s,b_1-r_1,\dots,b_m-r_m).
$$
\end{definition}

\begin{example}
Let us check if $\diag(6,6,6,3,1)$ is $4$-reducible and find its $4$-reduction.
We have $(b_1,b_2,b_3,b_4)=(6,6,3,1)$. Beginning with $r_4$ we can see
that $b_4=1<4$, so $r_4=b_4=1$, the same for $r_3=b_3=3$. Now $b_2\geq 4$, so
we take $\{1,2,3,4\} \setminus \{r_3,r_4\} = \{2,4\}$ and $r_2=4$, which is maximal.
The same applies for $r_1=2$. We can see that
$(r_1,r_2,r_3,r_4)=(2,4,3,1)$, hence $\diag(6,6,6,3,1)$ is $4$-reducible
and $\red_4(\diag(6,6,6,3,1))=\diag(6,4,2)$.
We present also another examples of reducing, and two diagrams which
are not reducible.
$$
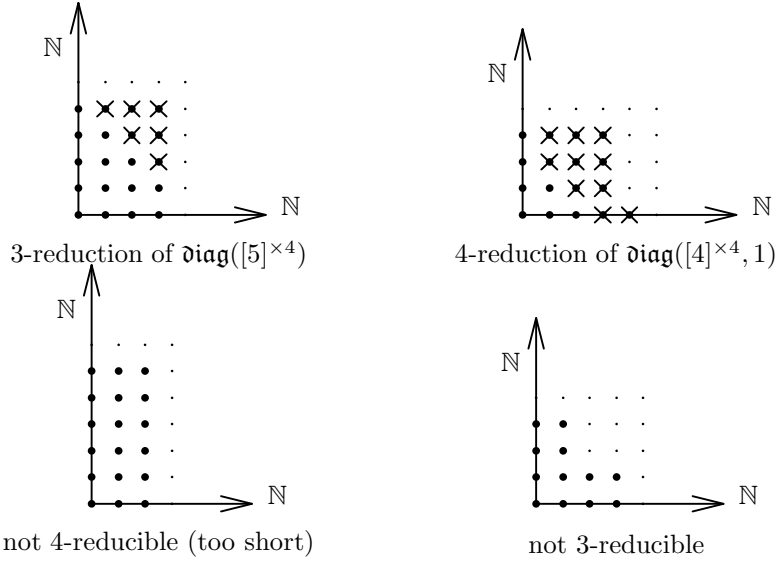
\begin{array}{ccc}

\begin{texdraw}
\drawdim pt
\move(0 0)
\fcir f:0 r:0.5
\move(0 10)
\fcir f:0 r:0.5
\move(0 20)
\fcir f:0 r:0.5
\move(0 30)
\fcir f:0 r:0.5
\move(0 40)
\fcir f:0 r:0.5
\move(0 50)
\fcir f:0 r:0.5
\move(10 0)
\fcir f:0 r:0.5
\move(10 10)
\fcir f:0 r:0.5
\move(10 20)
\fcir f:0 r:0.5
\move(10 30)
\fcir f:0 r:0.5
\move(10 40)
\fcir f:0 r:0.5
\move(10 50)
\fcir f:0 r:0.5
\move(20 0)
\fcir f:0 r:0.5
\move(20 10)
\fcir f:0 r:0.5
\move(20 20)
\fcir f:0 r:0.5
\move(20 30)
\fcir f:0 r:0.5
\move(20 40)
\fcir f:0 r:0.5
\move(20 50)
\fcir f:0 r:0.5
\move(30 0)
\fcir f:0 r:0.5
\move(30 10)
\fcir f:0 r:0.5
\move(30 20)
\fcir f:0 r:0.5
\move(30 30)
\fcir f:0 r:0.5
\move(30 40)
\fcir f:0 r:0.5
\move(30 50)
\fcir f:0 r:0.5
\move(40 0)
\fcir f:0 r:0.5
\move(40 10)
\fcir f:0 r:0.5
\move(40 20)
\fcir f:0 r:0.5
\move(40 30)
\fcir f:0 r:0.5
\move(40 40)
\fcir f:0 r:0.5
\move(40 50)
\fcir f:0 r:0.5
\arrowheadtype t:V
\move(0 0)
\avec(70 0)
\move(0 0)
\avec(0 80)
\htext(76 0){$\mathbb{N}$}
\htext(-13 60){$\mathbb{N}$}
\move(0 0)
\fcir f:0 r:1.5
\move(0 10)
\fcir f:0 r:1.5
\move(0 20)
\fcir f:0 r:1.5
\move(0 30)
\fcir f:0 r:1.5
\move(0 40)
\fcir f:0 r:1.5
\move(10 0)
\fcir f:0 r:1.5
\move(10 10)
\fcir f:0 r:1.5
\move(10 20)
\fcir f:0 r:1.5
\move(10 30)
\fcir f:0 r:1.5
\move(10 40)
\fcir f:0 r:1.5
\move(20 0)
\fcir f:0 r:1.5
\move(20 10)
\fcir f:0 r:1.5
\move(20 20)
\fcir f:0 r:1.5
\move(20 30)
\fcir f:0 r:1.5
\move(20 40)
\fcir f:0 r:1.5
\move(30 0)
\fcir f:0 r:1.5
\move(30 10)
\fcir f:0 r:1.5
\move(30 20)
\fcir f:0 r:1.5
\move(30 30)
\fcir f:0 r:1.5
\move(30 40)
\fcir f:0 r:1.5
\move(10 40)
\move(7 37)
\lvec(13 43)
\move(13 37)
\lvec(7 43)
\move(20 40)
\move(17 37)
\lvec(23 43)
\move(23 37)
\lvec(17 43)
\move(20 30)
\move(17 27)
\lvec(23 33)
\move(23 27)
\lvec(17 33)
\move(30 40)
\move(27 37)
\lvec(33 43)
\move(33 37)
\lvec(27 43)
\move(30 30)
\move(27 27)
\lvec(33 33)
\move(33 27)
\lvec(27 33)
\move(30 20)
\move(27 17)
\lvec(33 23)
\move(33 17)
\lvec(27 23)
\textref h:C v:C
\htext(30 -15){$3$-reduction of $\diag([5]^{\times 4})$}
\end{texdraw}

&
\hspace{2cm}
&

\begin{texdraw}
\drawdim pt
\move(0 0)
\fcir f:0 r:0.5
\move(0 10)
\fcir f:0 r:0.5
\move(0 20)
\fcir f:0 r:0.5
\move(0 30)
\fcir f:0 r:0.5
\move(0 40)
\fcir f:0 r:0.5
\move(10 0)
\fcir f:0 r:0.5
\move(10 10)
\fcir f:0 r:0.5
\move(10 20)
\fcir f:0 r:0.5
\move(10 30)
\fcir f:0 r:0.5
\move(10 40)
\fcir f:0 r:0.5
\move(20 0)
\fcir f:0 r:0.5
\move(20 10)
\fcir f:0 r:0.5
\move(20 20)
\fcir f:0 r:0.5
\move(20 30)
\fcir f:0 r:0.5
\move(20 40)
\fcir f:0 r:0.5
\move(30 0)
\fcir f:0 r:0.5
\move(30 10)
\fcir f:0 r:0.5
\move(30 20)
\fcir f:0 r:0.5
\move(30 30)
\fcir f:0 r:0.5
\move(30 40)
\fcir f:0 r:0.5
\move(40 0)
\fcir f:0 r:0.5
\move(40 10)
\fcir f:0 r:0.5
\move(40 20)
\fcir f:0 r:0.5
\move(40 30)
\fcir f:0 r:0.5
\move(40 40)
\fcir f:0 r:0.5
\move(50 0)
\fcir f:0 r:0.5
\move(50 10)
\fcir f:0 r:0.5
\move(50 20)
\fcir f:0 r:0.5
\move(50 30)
\fcir f:0 r:0.5
\move(50 40)
\fcir f:0 r:0.5
\arrowheadtype t:V
\move(0 0)
\avec(80 0)
\move(0 0)
\avec(0 70)
\htext(86 0){$\mathbb{N}$}
\htext(-13 50){$\mathbb{N}$}
\move(0 0)
\fcir f:0 r:1.5
\move(0 10)
\fcir f:0 r:1.5
\move(0 20)
\fcir f:0 r:1.5
\move(0 30)
\fcir f:0 r:1.5
\move(10 0)
\fcir f:0 r:1.5
\move(10 10)
\fcir f:0 r:1.5
\move(10 20)
\fcir f:0 r:1.5
\move(10 30)
\fcir f:0 r:1.5
\move(20 0)
\fcir f:0 r:1.5
\move(20 10)
\fcir f:0 r:1.5
\move(20 20)
\fcir f:0 r:1.5
\move(20 30)
\fcir f:0 r:1.5
\move(30 0)
\fcir f:0 r:1.5
\move(30 10)
\fcir f:0 r:1.5
\move(30 20)
\fcir f:0 r:1.5
\move(30 30)
\fcir f:0 r:1.5
\move(40 0)
\fcir f:0 r:1.5
\move(10 30)
\move(7 27)
\lvec(13 33)
\move(13 27)
\lvec(7 33)
\move(10 20)
\move(7 17)
\lvec(13 23)
\move(13 17)
\lvec(7 23)
\move(20 30)
\move(17 27)
\lvec(23 33)
\move(23 27)
\lvec(17 33)
\move(20 20)
\move(17 17)
\lvec(23 23)
\move(23 17)
\lvec(17 23)
\move(20 10)
\move(17 7)
\lvec(23 13)
\move(23 7)
\lvec(17 13)
\move(30 30)
\move(27 27)
\lvec(33 33)
\move(33 27)
\lvec(27 33)
\move(30 20)
\move(27 17)
\lvec(33 23)
\move(33 17)
\lvec(27 23)
\move(30 10)
\move(27 7)
\lvec(33 13)
\move(33 7)
\lvec(27 13)
\move(30 0)
\move(27 -3)
\lvec(33 3)
\move(33 -3)
\lvec(27 3)
\move(40 0)
\move(37 -3)
\lvec(43 3)
\move(43 -3)
\lvec(37 3)
\textref h:C v:C
\htext(35 -15){$4$-reduction of $\diag([4]^{\times 4},1)$}
\end{texdraw}

\\

\begin{texdraw}
\drawdim pt
\move(0 0)
\fcir f:0 r:0.5
\move(0 10)
\fcir f:0 r:0.5
\move(0 20)
\fcir f:0 r:0.5
\move(0 30)
\fcir f:0 r:0.5
\move(0 40)
\fcir f:0 r:0.5
\move(0 50)
\fcir f:0 r:0.5
\move(0 60)
\fcir f:0 r:0.5
\move(10 0)
\fcir f:0 r:0.5
\move(10 10)
\fcir f:0 r:0.5
\move(10 20)
\fcir f:0 r:0.5
\move(10 30)
\fcir f:0 r:0.5
\move(10 40)
\fcir f:0 r:0.5
\move(10 50)
\fcir f:0 r:0.5
\move(10 60)
\fcir f:0 r:0.5
\move(20 0)
\fcir f:0 r:0.5
\move(20 10)
\fcir f:0 r:0.5
\move(20 20)
\fcir f:0 r:0.5
\move(20 30)
\fcir f:0 r:0.5
\move(20 40)
\fcir f:0 r:0.5
\move(20 50)
\fcir f:0 r:0.5
\move(20 60)
\fcir f:0 r:0.5
\move(30 0)
\fcir f:0 r:0.5
\move(30 10)
\fcir f:0 r:0.5
\move(30 20)
\fcir f:0 r:0.5
\move(30 30)
\fcir f:0 r:0.5
\move(30 40)
\fcir f:0 r:0.5
\move(30 50)
\fcir f:0 r:0.5
\move(30 60)
\fcir f:0 r:0.5
\arrowheadtype t:V
\move(0 0)
\avec(60 0)
\move(0 0)
\avec(0 90)
\htext(66 0){$\mathbb{N}$}
\htext(-13 70){$\mathbb{N}$}
\move(0 0)
\fcir f:0 r:1.5
\move(0 10)
\fcir f:0 r:1.5
\move(0 20)
\fcir f:0 r:1.5
\move(0 30)
\fcir f:0 r:1.5
\move(0 40)
\fcir f:0 r:1.5
\move(0 50)
\fcir f:0 r:1.5
\move(10 0)
\fcir f:0 r:1.5
\move(10 10)
\fcir f:0 r:1.5
\move(10 20)
\fcir f:0 r:1.5
\move(10 30)
\fcir f:0 r:1.5
\move(10 40)
\fcir f:0 r:1.5
\move(10 50)
\fcir f:0 r:1.5
\move(20 0)
\fcir f:0 r:1.5
\move(20 10)
\fcir f:0 r:1.5
\move(20 20)
\fcir f:0 r:1.5
\move(20 30)
\fcir f:0 r:1.5
\move(20 40)
\fcir f:0 r:1.5
\move(20 50)
\fcir f:0 r:1.5
\textref h:C v:C
\htext(25 -15){not $4$-reducible (too short)}
\end{texdraw}

&
\hspace{2cm}
&

\begin{texdraw}
\drawdim pt
\move(0 0)
\fcir f:0 r:0.5
\move(0 10)
\fcir f:0 r:0.5
\move(0 20)
\fcir f:0 r:0.5
\move(0 30)
\fcir f:0 r:0.5
\move(0 40)
\fcir f:0 r:0.5
\move(10 0)
\fcir f:0 r:0.5
\move(10 10)
\fcir f:0 r:0.5
\move(10 20)
\fcir f:0 r:0.5
\move(10 30)
\fcir f:0 r:0.5
\move(10 40)
\fcir f:0 r:0.5
\move(20 0)
\fcir f:0 r:0.5
\move(20 10)
\fcir f:0 r:0.5
\move(20 20)
\fcir f:0 r:0.5
\move(20 30)
\fcir f:0 r:0.5
\move(20 40)
\fcir f:0 r:0.5
\move(30 0)
\fcir f:0 r:0.5
\move(30 10)
\fcir f:0 r:0.5
\move(30 20)
\fcir f:0 r:0.5
\move(30 30)
\fcir f:0 r:0.5
\move(30 40)
\fcir f:0 r:0.5
\move(40 0)
\fcir f:0 r:0.5
\move(40 10)
\fcir f:0 r:0.5
\move(40 20)
\fcir f:0 r:0.5
\move(40 30)
\fcir f:0 r:0.5
\move(40 40)
\fcir f:0 r:0.5
\arrowheadtype t:V
\move(0 0)
\avec(70 0)
\move(0 0)
\avec(0 70)
\htext(76 0){$\mathbb{N}$}
\htext(-13 50){$\mathbb{N}$}
\move(0 0)
\fcir f:0 r:1.5
\move(0 10)
\fcir f:0 r:1.5
\move(0 20)
\fcir f:0 r:1.5
\move(0 30)
\fcir f:0 r:1.5
\move(10 0)
\fcir f:0 r:1.5
\move(10 10)
\fcir f:0 r:1.5
\move(10 20)
\fcir f:0 r:1.5
\move(10 30)
\fcir f:0 r:1.5
\move(20 0)
\fcir f:0 r:1.5
\move(20 10)
\fcir f:0 r:1.5
\move(30 0)
\fcir f:0 r:1.5
\move(30 10)
\fcir f:0 r:1.5
\textref h:C v:C
\htext(30 -15){not $3$-reducible}
\end{texdraw}

\end{array}
$$
\end{example}

\begin{definition}
Let $D \subset \Z^2$ be a finite set, let $m_1,\dots,m_r$ be non-negative
integers. We will identify points $(\alpha,\beta) \in \Z^2$ with
monomials $x^\alpha y^\beta \in \field[x,x^{-1},y,y^{-1}]$.
Take points $p_1,\dots,p_r$ in general position in $\field^2$
and define the vector space (over $\field$)
\begin{multline*}
\sys(D;m_1,\dots,m_r) = \\
\{ f \in \field[x,x^{-1},y,y^{-1}] : \supp(f) \subset D, \,
\mult_{p_j}f \geq m_j \text{ for } j=1,\dots,r\}.
\end{multline*}
Put
\begin{align*}
\dim \sys(D;m_1,\dots,m_r) & = \dim_{\field} \sys(D;m_1,\dots,m_r) - 1, \\
\vdim \sys(D;m_1,\dots,m_r) & = \#D - 1 - \sum_{j=1}^{r} \binom{m_j+1}{2}, \\
\edim \sys(D;m_1,\dots,m_r) & = \max \{ \vdim \sys(D;m_1,\dots,m_r), -1 \}.
\end{align*}
We say that $\sys(D;m_1,\dots,m_r)$ is \emph{special} if 
$$\dim\sys(D;m_1,\dots,m_r) > \edim\sys(D;m_1,\dots,m_r).$$
\end{definition}

\begin{lemma}
\label{diagshift}
Let $D \subset \Z^2$ be finite, let $\varphi$ be one of
the following maps
\begin{gather*}
\Z^2 \ni (a,b) \mapsto (b,a) \in \Z^2, \qquad
\Z^2 \ni (a,b) \mapsto (a,-b) \in \Z^2, \\
\Z^2 \ni (a,b) \mapsto (a,b+a) \in \Z^2, \qquad
\Z^2 \ni (a,b) \mapsto (a,b+c) \in \Z^2,
\end{gather*}
where $c \in \Z$.
Then
$$
\dim \sys(D;m_1,\dots,m_r) = \dim \sys(\varphi(D);m_1,\dots,m_r)
$$
for any non-negative integers $m_1,\dots,m_r$.
\end{lemma}

\begin{proof}
Pick generic points $p_1,\dots,p_r$. We can assume that
$(p_j)_x \neq 0, (p_j)_y \neq 0$, where by $(p)_x$, $(p)_y$ we denote
the first and the second coordinate of a point $p$, respectively.
It can be shown that the following linear maps 
from $\field[x,x^{-1},y,y^{-1}] \to \field[x,x^{-1},y,y^{-1}]$
given by
\begin{gather*}
x^ay^b \mapsto y^ax^b , \qquad
x^ay^b \mapsto x^a\frac{1}{y^b} ,\\
x^ay^b \mapsto (xy)^ay^b, \qquad
x^ay^b \mapsto y^cx^ay^b
\end{gather*}
induce isomorphisms of $\sys(D;m_1,\dots,m_r)$ with
$\sys(\varphi(D);m_1,\dots,m_r)$, where in the last systems the coordinates
of base points $q_1,\dots,q_r$ are given by
\begin{gather*}
(q_j)_x = (p_j)_y, (q_j)_y = (p_j)_x, \qquad
(q_j)_x = (p_j)_x, (q_j)_y = \frac{1}{(p_j)_y}, \\
(q_j)_x = \frac{(p_j)_x}{(p_j)_y}, (q_j)_y = (p_j)_y, \qquad
(q_j)_x = (p_j)_x, (q_j)_y = (p_j)_y,\end{gather*}
respectively.
To see this, observe that $\frac{\partial}{\partial y^n}(f(x,y^{-1}))$
can be written as a linear combination (over $\field[y]$)
of $(\frac{\partial f}{\partial y^j})(x,y^{-1})$ for $j=0,\dots,n$.
This suffices to complete the case $(a,b) \mapsto (a,-b)$.
The rest of the proof is left to the reader.
\end{proof}

\begin{proposition}
\label{diagforhi}
Let $n$, $b$, $a$, $m_1,\dots,m_r$ be non-negative integers.
Let
$$D = \diag([1]^{\times n},[2]^{\times n},\dots,[b]^{\times n},[b+1]^{\times (a+1)}).$$
Then
\begin{align*}
\dim \sys_{n}(a,b;m_1,\dots,m_r) & = \dim \sys(D;m_1,\dots,m_r),\\
\vdim \sys_{n}(a,b;m_1,\dots,m_r) & = \vdim \sys(D;m_1,\dots,m_r).
\end{align*}
\end{proposition}

\begin{proof}
The Hirzebruch surface $\FF_n$ is a toric surface given by
the fan generated by $v_1=e_1$, $v_2=e_2$, $v_3=-e_1+ne_2$, $v_4=-e_2$
(Figure \ref{genforhi}; more theory on toric varieties can be found in \cite{Fulton}).
\begin{figure}[ht!]
\begin{texdraw}
\drawdim pt
\arrowheadtype t:V
\move(0 0)
\avec(60 0)
\move(0 0)
\avec(0 80)
\htext(66 0){$\mathbb{N}$}
\htext(-13 65){$\mathbb{N}$}
\linewd 1
\move(0 0)
\avec(20 0)
\move(0 0)
\avec(0 20)
\move(0 0)
\avec(-20 60)
\move(0 0)
\avec(0 -20)
\htext(20 5){{\scriptsize $v_1$}}
\htext(5 20){{\scriptsize $v_2$}}
\htext(-25 60){{\scriptsize $v_3$}}
\htext(5 -20){{\scriptsize $v_4$}}
\end{texdraw}
\caption{Toric fan for $\FF_n$}\label{genforhi}
\end{figure}
The class $F$ is given by $v_1$ as the class of the curve corresponding
to $v_1$, similarly $H$ is given by $v_4$.
Now, for the line bundle $aF+bH$ we have
$$\text{global sections of } \, aF+bH = \bigoplus_{u \in P_{aF+bH}} \field \chi^{u},$$
where
$$
P_{aF+bH} = \{ u \in \Z^2 : \langle u,v_i \rangle \geq -a_i \} 
= \{ (u_x,u_y) \in \Z^2 : -a\leq u_x \leq nu_y, 0 \leq u_y \leq b \}
$$
(see Figure \ref{fanshape}),
\begin{figure}[ht!]
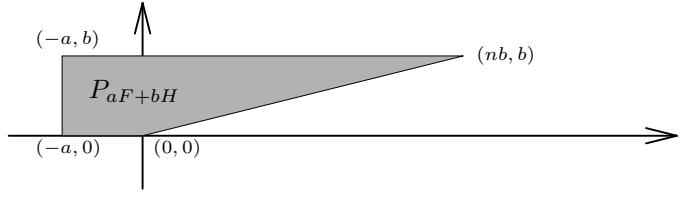

\begin{texdraw}
\drawdim pt
\arrowheadtype t:V
\move(-50 0)
\avec(200 0)
\move(0 -20)
\avec(0 50)
\linewd 0.1
\move(-30 0)
\lvec(0 0)
\lvec(120 30)
\lvec(-30 30)
\lvec(-30 0)
\lfill f:0.7
\htext(4 -8){{\scriptsize $(0,0)$}}
\htext(125 27){{\scriptsize $(nb,b)$}}
\htext(-40 33){{\scriptsize $(-a,b)$}}
\htext(-40 -8){{\scriptsize $(-a,0)$}}
\htext(-20 12){$P_{aF+bH}$}
\end{texdraw}
\caption{Global sections of $aF+bH$}\label{fanshape}
\end{figure}
since
$a_1=a$, $a_4=b$, $a_2=a_3=0$.
The element $\chi^{(u_x,u_y)}$ is identified with the monomial $x^{u_x}y^{u_y}$ 
in $\field[x,x^{-1},y,y^{-1}]$.
Since there is a Zariski open, non-empty set $U$ on $\FF_n$
such that $U$ is affine with coordinate ring
$\field[x,y,x^{-1},y^{-1}]$, taking $D=P_{aF+bH}$ and points in general
position we have
$$\dim \sys_n(a,b;m_1,\dots,m_r) = \dim \sys(D;m_1,\dots,m_r).$$
By Lemma \ref{diagshift}, 
applying $\varphi: (p,q) \mapsto (p+a+n(b-q),b-q)$ we can transform
$P_{aF+bH}$ into the diagram contained in the trapezoid with vertices
$(0,0)$, $(a+nb,0)$, $(a+nb,b)$, $(nb,b)$ (see Figure \ref{fanshapeD}),
\begin{figure}[ht!]
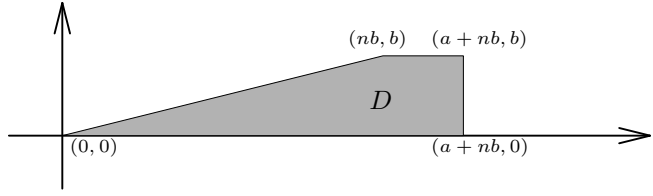

\begin{texdraw}
\drawdim pt
\arrowheadtype t:V
\move(-20 0)
\avec(220 0)
\move(0 -20)
\avec(0 50)
\linewd 0.1
\move(0 0)
\lvec(150 0)
\lvec(150 30)
\lvec(120 30)
\lvec(0 0)
\lfill f:0.7
\htext(3 -8){{\scriptsize $(0,0)$}}
\htext(138 -8){{\scriptsize $(a+nb,0)$}}
\htext(138 33){{\scriptsize $(a+nb,b)$}}
\htext(107 33){{\scriptsize $(nb,b)$}}
\htext(115 10){$D$}
\end{texdraw}
\caption{The image of $P_{aF+bH}$ by $\varphi$}\label{fanshapeD}
\end{figure}
which completes the proof.
For virtual dimension we use the fact that $\dim \sys_n(a,b)= \# D - 1$.
\end{proof}

Before formulating the next proposition, we need one additional notation.
This notation will be used only in the following Proposition \ref{diagforp1p1}
and then
in Proposition \ref{setpb}.

\begin{definition}
Let $a_1,\dots,a_s$, $u_1,\dots,u_s$ be non-negative integers.
Define
$$\diag(a_1^{\up u_1},\dots,a_s^{\up u_s}) = \{ (x,y) \in \N^2 : u_{x+1} \leq y < a_{x+1}+u_{x+1} \}.$$
\end{definition}

\begin{proposition}
\label{diagforp1p1}
Let $b$, $a$, $m_1,\dots,m_r$ be non-negative integers.
Let
$$D = \diag(1,\dots,b,[b+1]^{\times (a-b+1)},b^{\up 1},(b-1)^{\up 2},\dots,1^{\up b}).$$
Then
\begin{align*}
\dim \sys_{0}(a,b;m_1,\dots,m_r) & = \dim \sys(D;m_1,\dots,m_r),\\
\vdim \sys_{0}(a,b;m_1,\dots,m_r) & = \vdim \sys(D;m_1,\dots,m_r).
\end{align*}
\end{proposition}

\begin{proof}
Let $G = \diag([a+1]^{\times (b+1)})$. From the previous
proof and Lemma \ref{diagshift} we have
$$\dim \sys(G;m_1,\dots,m_r) = \dim \sys_{0}(a,b;m_1,\dots,m_r).$$
By applying
$\varphi : (p,q) \mapsto (p+q,q)$ we obtain $D=\varphi(G)$ (see Figure \ref{p1p1shape}).
\begin{figure}[ht!]
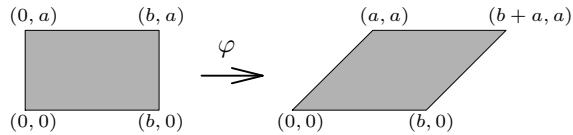

\begin{texdraw}
\drawdim pt
\arrowheadtype t:V

\move(66 13)
\avec(89 13)
\htext(72 20){$\varphi$}

\linewd 0.1
\move(0 0)
\lvec(50 0)
\lvec(50 30)
\lvec(0 30)
\lvec(0 0)
\lfill f:0.7
\htext(-6 -8){{\scriptsize $(0,0)$}}
\htext(42 -8){{\scriptsize $(b,0)$}}
\htext(42 32){{\scriptsize $(b,a)$}}
\htext(-6 32){{\scriptsize $(0,a)$}}

\move(100 0)
\lvec(150 0)
\lvec(180 30)
\lvec(130 30)
\lvec(100 0)
\lfill f:0.7
\htext(94 -8){{\scriptsize $(0,0)$}}
\htext(143 -8){{\scriptsize $(b,0)$}}
\htext(173 32){{\scriptsize $(b+a,a)$}}
\htext(125 32){{\scriptsize $(a,a)$}}

\end{texdraw}
\caption{Global sections of $aF+bH$ for $\FF_0$}\label{p1p1shape}
\end{figure}
\end{proof}

\begin{theorem}
\label{Mred}
Let $D$, $G$ be diagrams, let $m\geq 1$, $p,q \geq 0$. If
\begin{itemize}
\item
the diagram $G$ can be obtained from $D$ by a sequence of $p$ $m$-reductions,
\item
the system $\sys(G;m^{\times q},m_1,\dots,m_r)$ is non-special
\end{itemize}
then the system $\sys(D;m^{\times (q+p)},m_1,\dots,m_r)$ is non-special.
\end{theorem}

\begin{proof}
The proof can be found in \cite[proof of Theorem 7 and of Proposition 18]{redmd}.
\end{proof}

\begin{example}
We will show that $\sys_2(2,3;3^{\times 4})$ is non-special.
By Proposition \ref{diagforhi} we have to show that
$\sys(D;3^{\times 4})$ is empty for $D=\diag(1,1,2,2,3,3,4,4,4)$.
By Theorem \ref{Mred} it is enough to show emptiness of
$\sys(G;3^{\times 2})$ for
$$G = \red_3(\red_3(D)) = \diag(1,1,2,2,3,3).$$
The last diagram is $3$-reducible, but then the system
$\sys(\diag(1,1,2,1,1);3)$ is non-empty. Instead, using Lemma \ref{diagshift},
we can change $\diag(1,1,2,2,3,3)$ into $\diag(2,4,6)$, which can be $3$-reduced
twice to the empty diagram.
\end{example}

\section{Speciality of some homogeneous systems}

\begin{proposition}
\label{easyf1}
The system $\sys_{1}(a,b;m^{\times r})$ for $m \leq 10$ is special if and
only if it is $-1$-special.
\end{proposition}

\begin{proof}
The surface $\mathbb F_1$ is the blow-up of $\mathbb P^2$ in one point,
so the proposition follows from \cite[Theorem 3]{Dumquasi}, where it is shown that
the SHGH Conjecture holds for quasi-homogeneous systems with
homogeneous multiplicity bounded by $10$.
\end{proof}

\begin{proposition}
\label{easyfn}
The system $\sys_{n}(a,b;m^{\times r})$ for $n \geq 2$ and $b \leq m+1 \leq 11$
is special if and only if it is $-1$-special.
\end{proposition}

\begin{proof}
The proof for the case
$r \leq n+1$ and arbitrary $b$ and $m$ can be found in \cite[Proposition 2.7]{Laf},
but there is serious mistake --- the line bundle $-K_S$ fails to be nef,
where $K_S$ is the canonical
bundle on the blow-up of $\FF_n$ at $r$ points. This is due to the fact
that $K_{\FF_n}^2=8$, so $(-K_S)^2 = 8-r$ and this number is negative
for $r > 8$ points.

We will give a proof different from that in \cite{Laf}.
Moreover, we will also consider the case $r \leq n+1$ separately and prove that
each system of this type is either non-special or $-1$-special without
our additional assumption that $m \leq 10$.

First of all, due to Proposition \ref{fromhitop2}, we will work with the planar system
$$\sys(a+(n+1)b;a+nb,m^{\times r}, \overline{b^{\times (n+1)}}).$$
During the proof we will write $(k)_{\geq 0}$ for $\max\{k,0\}$.

The first case to consider is $r > n+1$.
By Cremona transformations
based on the point with the greatest multiplicity, points lying on a line
and points with multiplicity $b$ we obtain
$$\sys(a+(n+1)b-(n+1)m;(a+nb-(n+1)m)_{\geq 0},m^{\times (r-n-1)},\overline{(b-m)_{\geq 0}^{\times (n+1)}}).$$
For $b-m \leq 0$ we are done by \cite{Dumquasi}.
If $b-m = 1$ then either the line $\ell$ supporting $n+1$ points is in the base locus,
hence the system is non-special or $-1$-special,
or each point lying on this line imposes an independent condition and
\begin{multline*}
\dim \sys(a+(n+1)(b-m);(a+nb-(n+1)m)_{\geq 0},m^{\times (r-n-1)},\overline{1^{\times (n+1)}}) \\
= \dim \sys(a+(n+1)(b-m);(a+nb-(n+1)m)_{\geq 0},m^{\times (r-n-1)},1^{\times (n+1)}).
\end{multline*}
If the last system is special then
$$
\sys(a+(n+1)(b-m);(a+nb-(n+1)m)_{\geq 0},m^{\times (r-n-1)})
$$
is special and then, by \cite{Dumquasi}, $-1$-special. 
This finishes the first case.

The second case is $r \leq n+1$. As above, by
Cremona we can consider
$$\sys(a+(n+1)b-rm;(a+nb-rm)_{\geq 0},\overline{(b-m)_{\geq 0}^{\times r},b^{\times (n+1-r)}}).$$
Now we can remove the fixed part (again consisting of lines) to obtain
either the empty system, or the system $L=\sys(d;m_0,\overline{m_1,\dots,m_s})$ of the same dimension,
where non-negative
integers $d,m_0,m_1,\dots,m_s$ satisfy
$$d \geq \sum_{j=1}^{s} m_j \quad \text{and} \quad d \geq m_0 + m_j \text{ for } j=1,\dots,s.$$
Moreover, we will assume that
$$m_1 \geq m_2 \geq \dots \geq m_s.$$
We will show that the last system is non-special, 
so by Proposition \ref{cremona} the system we begin with, is either $-1$-special
or non-special.

By a suitable projective change of coordinates we can assume that
$p_0=(0:1:0)$ and collinear points $q_1,\dots,q_s$ have coordinates
$(w_1:0:1),\dots,(w_s:0:1)$.
Take $C \in L$ defined by a polynomial $f$. Then $f$ is generated 
by monomials contained in the set
\begin{align*}
D & = \{x^ay^bz^c : a+b+c=d\} = U \cup G \cup (D \setminus (U \cup G)), \\
\intertext{where}
U & = \{x^ay^bz^c : a+b+c=d, \, b > d-m\}, \\
G & = \big\{x^ay^bz^c : a+b+c=d, \, a < \sum_{j} (m_j-b)_{\geq 0} \big\}
\end{align*}
(see Figure \ref{figshade}; the picture is drawn after dehomogenizing with
respect to $z$).
By our assumption on multiplicities we know that
$$\# U = \binom{m_0+1}{2}, \quad \# G = \sum_{j=1}^{s} \binom{m_j+1}{2}, \quad
U \cap G = \varnothing.$$
\begin{figure}[ht!]
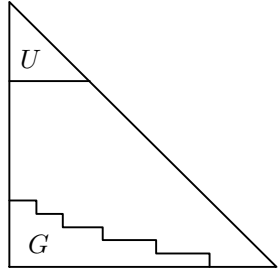

\begin{texdraw}
\drawdim pt
\move(0 0)
\lvec(100 0)
\lvec(0 100)
\lvec(0 0)
\move(0 70)
\lvec(30 70)
\htext(4 75){$U$}
\move(0 25)
\lvec(10 25)
\lvec(10 20)
\lvec(20 20)
\lvec(20 15)
\lvec(35 15)
\lvec(35 10)
\lvec(55 10)
\lvec(55 5)
\lvec(75 5)
\lvec(75 0)
\htext(7 5){$G$}
\end{texdraw}
\caption{Division of $D$ into $U$ and $G$}\label{figshade}
\end{figure}
It is enough to show that for fixed coefficients standing by monomials
from $D \setminus (U \cup G)$ there exists exactly one $f$, which defines a curve in $L$.
Indeed, we can see that $\supp(f) \cap U = \varnothing$. Let
$f=f_1+f_2+f_3$, where $\supp(f_1) \in G \cap \{x^ay^bz^c : b=0\}$,
so $\# \supp(f_1) \leq \sum_j m_j$,
$\supp(f_2) \in (D \setminus G) \cap \{a^ay^bz^c : b=0\}$,
$\supp(f_3) \in D \cap \{x^ay^bz^c : b > 0\}$.
We have
\begin{align*}
0 & = \frac{\partial^jf}{\partial x^j}(p_k)=
\frac{\partial^jf_1}{\partial x^j}(w_k:0:1)+
\frac{\partial^jf_2}{\partial x^j}(w_k:0:1)+
\frac{\partial^jf_3}{\partial x^j}(w_k:0:1) \\
& = \frac{\partial^j(f_1(x:0:1))}{\partial x^j}(w_k)+
\frac{\partial^jf_2}{\partial x^j}(w_k:0:1)
\end{align*}
for $j=0,\dots,m_k-1$, $k = 1,\dots,s$.
Since $\frac{\partial^jf_2}{\partial x^j}(w_k:0:1)$ is uniquely determined,
we use one dimensional interpolation
to uniquely determine $f_1$.
For other monomials in $G$ we deduce in the analogous way, using induction.
Namely, we assume that the coefficients standing by monomials
$$\{x^ay^bz^c : x^ay^bz^c \in G, \, b \leq n-1\}$$
are uniquely determined, and we will show the same for
$$\{x^ay^bz^c : x^ay^bz^c \in G, \, b = n\}.$$
Again, let $f=f_1+f_2+f_3+f_4$ for
$\supp(f_1) \in G \cap \{x^ay^bz^c : b=n\}$,
$\supp(f_2) \in (D \setminus G) \cap \{x^ay^bz^c : b=n\}$,
$\supp(f_3) \in \{x^ay^bz^c : b < n\}$,
$\supp(f_4) \in \{x^ay^bz^c : b > n\}$.
We have
\begin{align*}
0 & =
\bigg(\frac{\partial^{j+n}f_1}{\partial x^j \partial y^n}+
\frac{\partial^{j+n}f_2}{\partial x^j \partial y^n}+
\frac{\partial^{j+n}f_3}{\partial x^j \partial y^n}+
\frac{\partial^{j+n}f_4}{\partial x^j \partial y^n}\bigg)(w_k:0:1) \\
& = n! \frac{\partial^j(f_1(x:0:1))}{\partial x^j}(w_k)+
\frac{\partial^{j+n}f_2}{\partial x^j \partial y^n}(w_k:0:1)+
\frac{\partial^{j+n}f_4}{\partial x^j \partial y^n}(w_k:0:1)
\end{align*}
for $k = 1,\dots,s$ such that $m_k > n$ and $j=0,\dots,m_k-n-1$.
By assumptions, the contribution from $f_2$ and $f_4$ is fixed and we use
interpolation.
\end{proof}

\begin{proposition}
\label{easyf0}
The system $\sys_{0}(a,b;m^{\times r})$ for $\min\{a,b\} \leq m \leq 10$
is special if and only if it is $-1$-special.
\end{proposition}

\begin{proof}
The above system is equivalent to $\sys(a+b;a,b,m^{\times r})$,
so, by Cremona, we get
$$\sys(a+b-m;(a-m)_{\geq 0},(b-m)_{\geq 0},m^{\times (r-1)})$$
and we use \cite{Dumquasi} (for quasihomogeneous) or \cite[Theorem 32]{md}
(for homogeneous system).
\end{proof}

\section{Sequences of reductions}

\begin{definition}
Let $m \geq 2$, let $h > m$ be integers. We say that $(a_1,\dots,a_{m-1})$ is
an \emph{admissible $h$-$(b_1,\dots,b_{m-1})$-tail for multiplicity $m$}
if there exists $k \geq 0$ such that the diagram
$D_1=\diag(a_1,\dots,a_{m-1})$ can be obtained
from $D_2=\diag([h]^{\times k},b_1,\dots,b_m)$ by a sequence of
$m$-reductions.
\end{definition}

\begin{example}
Let us show how one can enumerate all admissible $8$-$([0]^{\times 3})$-tails for
multiplicity $4$. We begin with $(a_1,a_2,a_3)$ equal to
$(0,0,0)$, $(8,0,0)$, $(8,8,0)$ and $(8,8,8)$, which can be obtained without
performing any reduction. Now observe that $\diag(8,8,8,8)$
can be $4$-reduced to $\diag(6,4,2)$, so $(6,4,2)$ is also an
admissible tail. Now, each time we have an admissible tail, we can
add $8$ at the beginning and reduce until the fourth number disappears.
This gives the sequence
\begin{gather*}
\diag(6,4,2) \addarr \diag(8,6,4,2) \redarr \diag(7,3) \addarr \diag(8,7,3) \addarr \diag(8,8,7,3) \\
\redarr \diag(7,6,3) \addarr \diag(8,7,6,3) \redarr \diag(7,5,2) \addarr (8,7,5,2)  \\
\redarr \diag(7,4,1) \addarr \diag(8,7,4,1) \redarr \diag(6,4) \addarr \diag(8,6,4) \\
\addarr \diag(8,8,6,4) \redarr \diag(7,6,3) \addarr \dots
\end{gather*}
Observe that now our procedure will follow the loop, so nothing more
will appear.
Each diagram standing before one of the $\addarr$ arrows gives an admissible
tail.
\end{example}

\begin{definition}
Let $m \geq 2$, let $a_1,\dots,a_m$ be such that $\diag(a_1,\dots,a_m)$
is $m$-reducible. Let $\ell, k$ be nonnegative integers.
We say that 
\emph{$\diag(b_1,\dots,b_m,[\symbx]^{\times k})$ is a symbolic $m$-reduction
of $\diag(a_1,\dots,a_m,[\symbx]^{\times \ell})$} ($\symbx$ is just
the symbol without value) if and only if there exists
$c_1,\dots,c_\ell$ such that $m+1 \geq c_1 \geq \dots \geq c_\ell$
and $\diag(a_1,\dots,a_m,c_1,\dots,c_\ell)$ is $m$-reducible,
$$\red_m(\diag(a_1,\dots,a_m,c_1,\dots,c_\ell)) = \diag(b_1,\dots,b_m,d_1,\dots,d_\ell)$$
(with, possibly, some of $b_j$'s and $d_j$'s equal 0) satisfying
$$
k = \begin{cases}
0, & d_1=0, \\
\max\{j : d_j > 0\}, & d_1 > 0.
\end{cases}
$$
\end{definition}

\begin{example}
Let us enumerate all symbolic $3$-reductions for $D=\diag(6,6,7,\symbx,\symbx)$
(i.e. $\ell=2$).
This is just a straightforward calculation:
\begin{gather*}
\diag(6,6,7,4,4) \redarr \diag(6,6,6,2,1) \to \diag(6,6,6,\symbx,\symbx), \\
\diag(6,6,7,4,3) \redarr \diag(6,6,6,2) \to \diag(6,6,6,\symbx), \\
\diag(6,6,7,4,2) \redarr \diag(6,6,6,1) \to \diag(6,6,6,\symbx), \\
\diag(6,6,7,4,1) \redarr \diag(6,6,5,1) \to \diag(6,6,5,\symbx), \\
\diag(6,6,7,4) \redarr \diag(6,5,5,1) \to \diag(6,5,5,\symbx), \\
\diag(6,6,7,3,3) \redarr \diag(6,6,6,1) \to \diag(6,6,6,\symbx), \\
\diag(6,6,7,3,2) \redarr \diag(6,6,6), \quad
\diag(6,6,7,3,1) \redarr \diag(6,6,5), \\
\diag(6,6,7,3) \redarr \diag(6,5,5), \quad
\diag(6,6,7,2,1) \redarr \diag(6,6,4), \\
\diag(6,6,7,2) \redarr \diag(6,5,4), \quad
\diag(6,6,7,1) \redarr \diag(6,4,4), \\
\diag(6,6,7) \redarr \diag(5,4,4).
\end{gather*}
\end{example}

\begin{definition}
Let $m \geq 2$, let $a_1,\dots,a_m$ be such that $\diag(a_1,\dots,a_m)$
is $m$-reducible. We say that \emph{$D=\diag(b_1,\dots,b_{m-1})$
is an admissible $(a_1,\dots,a_m)$-tail for multiplicity $m$}
if $D$ can be obtained from $\diag(a_1,\dots,a_m,[\symbx]^{\times (m-1)})$
by a sequence of symbolic $m$-reductions.
\end{definition}

\begin{example}
Let us show that $\diag(6,4)$ is an admissible $(6,6,7)$-tail for
multiplicity $3$. Indeed,
\begin{gather*}
\diag(6,6,7,\symbx,\symbx) \symbredarr \diag(6,6,6,\symbx,\symbx)
\symbredarr \diag(6,6,5,\symbx,\symbx) \\
\symbredarr \diag(6,6,4,\symbx,\symbx) \symbredarr \diag(6,6,3,\symbx) \symbredarr \diag(6,4).
\end{gather*}
\end{example}

Observe that ``$\redarr$'' acts as a function, while ``$\symbredarr$''
is only a relation --- there are, usually, several possible
symbolic reductions. We will use $\symbred(D)$ to denote the set of all
symbolic reductions of $D$. Given $D$, we can find all elements in
$\symbred(D)$ in an algorithmic way, and since
every admissible $(a_1,\dots,a_m)$-tail belongs to
$$\symbred(\dots \symbred(\diag(a_1,\dots,a_m,[\symbx]^{\times (m-1)}))\dots),$$
where the number of symbolic reductions is bounded from above
(the bound depends on $m$ and $(a_1,\dots,a_m)$), we can enumerate
all admissible $(a_1,\dots,a_m)$-tails.

\begin{proposition}
\label{condonred}
Let $m \geq 2$, let $D=\diag(a_1,\dots,a_k)$ for 
$k \geq m$ and
\begin{itemize}
\item
$a_1 \geq m+1$,
\item
$a_{j+1} \in \{a_j-1,a_j,a_j+1\}$ for $j=1,\dots,k-1$,
\item
if $a_{j+1}=a_j+1$ then $a_j \geq 2m-1$ for $j=1,\dots,k-1$,
\item
if $a_{j+1}=a_j$ then $a_j \geq m$ for $j=1,\dots,k-1$.
\end{itemize}
Then $D$ can be reduced by a sequence of $m$-reductions
to a diagram $G$, which is an admissible $(a_1,\dots,a_m)$-tail
for multiplicity $m$.
\end{proposition}

\begin{proof}
We will reduce succesively, beginning with $D$. Let us
assume that we have obtained $\diag(b_1,\dots,b_k)$, which cannot
be $m$-reduced. Hence
$$0 < b_j \leq b_{j+1} \leq m-1$$
for some $j$.
We have three possibilities:
\begin{itemize}
\item
$a_{j+1} = a_j - 1$. Since the $(j+1)$th layer would be reduced stronger than
$j$th, we would have $b_{j+1} < b_j$;
\item
$a_{j+1} = a_j$. Now $b_j < a_j$ since $a_j \geq m$, so the $j$th layer must be
reduced at least once and again $b_{j+1} < b_j$;
\item
$a_{j+1} = a_j+1$. Now the difference between $a_j$ and $b_j$ is at least
$m$, so the $j$th layer must be reduced at least twice and again $b_{j+1} < b_j$.
\end{itemize}
So we can reduce $D$ to some $G=(c_1,\dots,c_{m-1})$. We have to show
that $G$ is an admissible $(a_1,\dots,a_m)$-tail for multiplicity $m$.

Define
$$\symb(\diag(d_1,\dots,d_m,d_{m+1},\dots,d_k)) = \diag(d_1,\dots,d_m,[\symbx]^{\times (k-m)}).$$
Observe that if
$E=\diag(d_1,\dots,d_m,e_1,\dots,e_k)$ is a diagram which can be $m$-reduced then
\begin{multline}
\label{eqsymb}
\symb(\red_m(E)) = \\
\symb(\red_m(\diag(d_1,\dots,d_m,\min\{e_1,m+1\},\dots,\min\{e_k,m+1\}))).
\end{multline}
The last diagram belongs to $\symbred_m(\symb(E))$.
Let $G=\red_m^{(p)}(D)$, where $\red_m^{(p)}$ denotes the $m$-reduction performed
$p$ times. From \eqref{eqsymb} we have
$$G=\symb(G) = \symb(\red_m^{(p)}(D)) \in \symbred^{(p)}(\symb(D)).$$
Since symbolic reductions performed on $\diag(a_1,\dots,a_m,[\symbx]^{\times s})$
for $s \geq m$ does not change $(a_1,\dots,a_m)$, we have
$$\symbred^{(p)}(\symb(D)) \subset \bigcup_{q=0}^{\infty} \symbred^{(q)}(\diag(a_1,\dots,a_m,[\symbx]^{\times (m-1)})),$$
which completes the proof.
\end{proof}

\section{Non-speciality by reductions}

\begin{figure}[ht!]
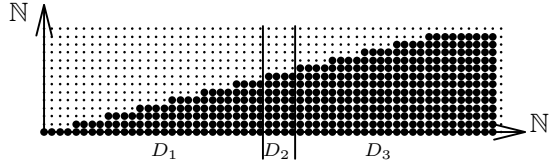

$$
\begin{texdraw}
\drawdim pt
\move(0 0)
\fcir f:0 r:0.5
\move(0 3)
\fcir f:0 r:0.5
\move(0 6)
\fcir f:0 r:0.5
\move(0 9)
\fcir f:0 r:0.5
\move(0 12)
\fcir f:0 r:0.5
\move(0 15)
\fcir f:0 r:0.5
\move(0 18)
\fcir f:0 r:0.5
\move(0 21)
\fcir f:0 r:0.5
\move(0 24)
\fcir f:0 r:0.5
\move(0 27)
\fcir f:0 r:0.5
\move(0 30)
\fcir f:0 r:0.5
\move(0 33)
\fcir f:0 r:0.5
\move(0 36)
\fcir f:0 r:0.5
\move(0 39)
\fcir f:0 r:0.5
\move(3 0)
\fcir f:0 r:0.5
\move(3 3)
\fcir f:0 r:0.5
\move(3 6)
\fcir f:0 r:0.5
\move(3 9)
\fcir f:0 r:0.5
\move(3 12)
\fcir f:0 r:0.5
\move(3 15)
\fcir f:0 r:0.5
\move(3 18)
\fcir f:0 r:0.5
\move(3 21)
\fcir f:0 r:0.5
\move(3 24)
\fcir f:0 r:0.5
\move(3 27)
\fcir f:0 r:0.5
\move(3 30)
\fcir f:0 r:0.5
\move(3 33)
\fcir f:0 r:0.5
\move(3 36)
\fcir f:0 r:0.5
\move(3 39)
\fcir f:0 r:0.5
\move(6 0)
\fcir f:0 r:0.5
\move(6 3)
\fcir f:0 r:0.5
\move(6 6)
\fcir f:0 r:0.5
\move(6 9)
\fcir f:0 r:0.5
\move(6 12)
\fcir f:0 r:0.5
\move(6 15)
\fcir f:0 r:0.5
\move(6 18)
\fcir f:0 r:0.5
\move(6 21)
\fcir f:0 r:0.5
\move(6 24)
\fcir f:0 r:0.5
\move(6 27)
\fcir f:0 r:0.5
\move(6 30)
\fcir f:0 r:0.5
\move(6 33)
\fcir f:0 r:0.5
\move(6 36)
\fcir f:0 r:0.5
\move(6 39)
\fcir f:0 r:0.5
\move(9 0)
\fcir f:0 r:0.5
\move(9 3)
\fcir f:0 r:0.5
\move(9 6)
\fcir f:0 r:0.5
\move(9 9)
\fcir f:0 r:0.5
\move(9 12)
\fcir f:0 r:0.5
\move(9 15)
\fcir f:0 r:0.5
\move(9 18)
\fcir f:0 r:0.5
\move(9 21)
\fcir f:0 r:0.5
\move(9 24)
\fcir f:0 r:0.5
\move(9 27)
\fcir f:0 r:0.5
\move(9 30)
\fcir f:0 r:0.5
\move(9 33)
\fcir f:0 r:0.5
\move(9 36)
\fcir f:0 r:0.5
\move(9 39)
\fcir f:0 r:0.5
\move(12 0)
\fcir f:0 r:0.5
\move(12 3)
\fcir f:0 r:0.5
\move(12 6)
\fcir f:0 r:0.5
\move(12 9)
\fcir f:0 r:0.5
\move(12 12)
\fcir f:0 r:0.5
\move(12 15)
\fcir f:0 r:0.5
\move(12 18)
\fcir f:0 r:0.5
\move(12 21)
\fcir f:0 r:0.5
\move(12 24)
\fcir f:0 r:0.5
\move(12 27)
\fcir f:0 r:0.5
\move(12 30)
\fcir f:0 r:0.5
\move(12 33)
\fcir f:0 r:0.5
\move(12 36)
\fcir f:0 r:0.5
\move(12 39)
\fcir f:0 r:0.5
\move(15 0)
\fcir f:0 r:0.5
\move(15 3)
\fcir f:0 r:0.5
\move(15 6)
\fcir f:0 r:0.5
\move(15 9)
\fcir f:0 r:0.5
\move(15 12)
\fcir f:0 r:0.5
\move(15 15)
\fcir f:0 r:0.5
\move(15 18)
\fcir f:0 r:0.5
\move(15 21)
\fcir f:0 r:0.5
\move(15 24)
\fcir f:0 r:0.5
\move(15 27)
\fcir f:0 r:0.5
\move(15 30)
\fcir f:0 r:0.5
\move(15 33)
\fcir f:0 r:0.5
\move(15 36)
\fcir f:0 r:0.5
\move(15 39)
\fcir f:0 r:0.5
\move(18 0)
\fcir f:0 r:0.5
\move(18 3)
\fcir f:0 r:0.5
\move(18 6)
\fcir f:0 r:0.5
\move(18 9)
\fcir f:0 r:0.5
\move(18 12)
\fcir f:0 r:0.5
\move(18 15)
\fcir f:0 r:0.5
\move(18 18)
\fcir f:0 r:0.5
\move(18 21)
\fcir f:0 r:0.5
\move(18 24)
\fcir f:0 r:0.5
\move(18 27)
\fcir f:0 r:0.5
\move(18 30)
\fcir f:0 r:0.5
\move(18 33)
\fcir f:0 r:0.5
\move(18 36)
\fcir f:0 r:0.5
\move(18 39)
\fcir f:0 r:0.5
\move(21 0)
\fcir f:0 r:0.5
\move(21 3)
\fcir f:0 r:0.5
\move(21 6)
\fcir f:0 r:0.5
\move(21 9)
\fcir f:0 r:0.5
\move(21 12)
\fcir f:0 r:0.5
\move(21 15)
\fcir f:0 r:0.5
\move(21 18)
\fcir f:0 r:0.5
\move(21 21)
\fcir f:0 r:0.5
\move(21 24)
\fcir f:0 r:0.5
\move(21 27)
\fcir f:0 r:0.5
\move(21 30)
\fcir f:0 r:0.5
\move(21 33)
\fcir f:0 r:0.5
\move(21 36)
\fcir f:0 r:0.5
\move(21 39)
\fcir f:0 r:0.5
\move(24 0)
\fcir f:0 r:0.5
\move(24 3)
\fcir f:0 r:0.5
\move(24 6)
\fcir f:0 r:0.5
\move(24 9)
\fcir f:0 r:0.5
\move(24 12)
\fcir f:0 r:0.5
\move(24 15)
\fcir f:0 r:0.5
\move(24 18)
\fcir f:0 r:0.5
\move(24 21)
\fcir f:0 r:0.5
\move(24 24)
\fcir f:0 r:0.5
\move(24 27)
\fcir f:0 r:0.5
\move(24 30)
\fcir f:0 r:0.5
\move(24 33)
\fcir f:0 r:0.5
\move(24 36)
\fcir f:0 r:0.5
\move(24 39)
\fcir f:0 r:0.5
\move(27 0)
\fcir f:0 r:0.5
\move(27 3)
\fcir f:0 r:0.5
\move(27 6)
\fcir f:0 r:0.5
\move(27 9)
\fcir f:0 r:0.5
\move(27 12)
\fcir f:0 r:0.5
\move(27 15)
\fcir f:0 r:0.5
\move(27 18)
\fcir f:0 r:0.5
\move(27 21)
\fcir f:0 r:0.5
\move(27 24)
\fcir f:0 r:0.5
\move(27 27)
\fcir f:0 r:0.5
\move(27 30)
\fcir f:0 r:0.5
\move(27 33)
\fcir f:0 r:0.5
\move(27 36)
\fcir f:0 r:0.5
\move(27 39)
\fcir f:0 r:0.5
\move(30 0)
\fcir f:0 r:0.5
\move(30 3)
\fcir f:0 r:0.5
\move(30 6)
\fcir f:0 r:0.5
\move(30 9)
\fcir f:0 r:0.5
\move(30 12)
\fcir f:0 r:0.5
\move(30 15)
\fcir f:0 r:0.5
\move(30 18)
\fcir f:0 r:0.5
\move(30 21)
\fcir f:0 r:0.5
\move(30 24)
\fcir f:0 r:0.5
\move(30 27)
\fcir f:0 r:0.5
\move(30 30)
\fcir f:0 r:0.5
\move(30 33)
\fcir f:0 r:0.5
\move(30 36)
\fcir f:0 r:0.5
\move(30 39)
\fcir f:0 r:0.5
\move(33 0)
\fcir f:0 r:0.5
\move(33 3)
\fcir f:0 r:0.5
\move(33 6)
\fcir f:0 r:0.5
\move(33 9)
\fcir f:0 r:0.5
\move(33 12)
\fcir f:0 r:0.5
\move(33 15)
\fcir f:0 r:0.5
\move(33 18)
\fcir f:0 r:0.5
\move(33 21)
\fcir f:0 r:0.5
\move(33 24)
\fcir f:0 r:0.5
\move(33 27)
\fcir f:0 r:0.5
\move(33 30)
\fcir f:0 r:0.5
\move(33 33)
\fcir f:0 r:0.5
\move(33 36)
\fcir f:0 r:0.5
\move(33 39)
\fcir f:0 r:0.5
\move(36 0)
\fcir f:0 r:0.5
\move(36 3)
\fcir f:0 r:0.5
\move(36 6)
\fcir f:0 r:0.5
\move(36 9)
\fcir f:0 r:0.5
\move(36 12)
\fcir f:0 r:0.5
\move(36 15)
\fcir f:0 r:0.5
\move(36 18)
\fcir f:0 r:0.5
\move(36 21)
\fcir f:0 r:0.5
\move(36 24)
\fcir f:0 r:0.5
\move(36 27)
\fcir f:0 r:0.5
\move(36 30)
\fcir f:0 r:0.5
\move(36 33)
\fcir f:0 r:0.5
\move(36 36)
\fcir f:0 r:0.5
\move(36 39)
\fcir f:0 r:0.5
\move(39 0)
\fcir f:0 r:0.5
\move(39 3)
\fcir f:0 r:0.5
\move(39 6)
\fcir f:0 r:0.5
\move(39 9)
\fcir f:0 r:0.5
\move(39 12)
\fcir f:0 r:0.5
\move(39 15)
\fcir f:0 r:0.5
\move(39 18)
\fcir f:0 r:0.5
\move(39 21)
\fcir f:0 r:0.5
\move(39 24)
\fcir f:0 r:0.5
\move(39 27)
\fcir f:0 r:0.5
\move(39 30)
\fcir f:0 r:0.5
\move(39 33)
\fcir f:0 r:0.5
\move(39 36)
\fcir f:0 r:0.5
\move(39 39)
\fcir f:0 r:0.5
\move(42 0)
\fcir f:0 r:0.5
\move(42 3)
\fcir f:0 r:0.5
\move(42 6)
\fcir f:0 r:0.5
\move(42 9)
\fcir f:0 r:0.5
\move(42 12)
\fcir f:0 r:0.5
\move(42 15)
\fcir f:0 r:0.5
\move(42 18)
\fcir f:0 r:0.5
\move(42 21)
\fcir f:0 r:0.5
\move(42 24)
\fcir f:0 r:0.5
\move(42 27)
\fcir f:0 r:0.5
\move(42 30)
\fcir f:0 r:0.5
\move(42 33)
\fcir f:0 r:0.5
\move(42 36)
\fcir f:0 r:0.5
\move(42 39)
\fcir f:0 r:0.5
\move(45 0)
\fcir f:0 r:0.5
\move(45 3)
\fcir f:0 r:0.5
\move(45 6)
\fcir f:0 r:0.5
\move(45 9)
\fcir f:0 r:0.5
\move(45 12)
\fcir f:0 r:0.5
\move(45 15)
\fcir f:0 r:0.5
\move(45 18)
\fcir f:0 r:0.5
\move(45 21)
\fcir f:0 r:0.5
\move(45 24)
\fcir f:0 r:0.5
\move(45 27)
\fcir f:0 r:0.5
\move(45 30)
\fcir f:0 r:0.5
\move(45 33)
\fcir f:0 r:0.5
\move(45 36)
\fcir f:0 r:0.5
\move(45 39)
\fcir f:0 r:0.5
\move(48 0)
\fcir f:0 r:0.5
\move(48 3)
\fcir f:0 r:0.5
\move(48 6)
\fcir f:0 r:0.5
\move(48 9)
\fcir f:0 r:0.5
\move(48 12)
\fcir f:0 r:0.5
\move(48 15)
\fcir f:0 r:0.5
\move(48 18)
\fcir f:0 r:0.5
\move(48 21)
\fcir f:0 r:0.5
\move(48 24)
\fcir f:0 r:0.5
\move(48 27)
\fcir f:0 r:0.5
\move(48 30)
\fcir f:0 r:0.5
\move(48 33)
\fcir f:0 r:0.5
\move(48 36)
\fcir f:0 r:0.5
\move(48 39)
\fcir f:0 r:0.5
\move(51 0)
\fcir f:0 r:0.5
\move(51 3)
\fcir f:0 r:0.5
\move(51 6)
\fcir f:0 r:0.5
\move(51 9)
\fcir f:0 r:0.5
\move(51 12)
\fcir f:0 r:0.5
\move(51 15)
\fcir f:0 r:0.5
\move(51 18)
\fcir f:0 r:0.5
\move(51 21)
\fcir f:0 r:0.5
\move(51 24)
\fcir f:0 r:0.5
\move(51 27)
\fcir f:0 r:0.5
\move(51 30)
\fcir f:0 r:0.5
\move(51 33)
\fcir f:0 r:0.5
\move(51 36)
\fcir f:0 r:0.5
\move(51 39)
\fcir f:0 r:0.5
\move(54 0)
\fcir f:0 r:0.5
\move(54 3)
\fcir f:0 r:0.5
\move(54 6)
\fcir f:0 r:0.5
\move(54 9)
\fcir f:0 r:0.5
\move(54 12)
\fcir f:0 r:0.5
\move(54 15)
\fcir f:0 r:0.5
\move(54 18)
\fcir f:0 r:0.5
\move(54 21)
\fcir f:0 r:0.5
\move(54 24)
\fcir f:0 r:0.5
\move(54 27)
\fcir f:0 r:0.5
\move(54 30)
\fcir f:0 r:0.5
\move(54 33)
\fcir f:0 r:0.5
\move(54 36)
\fcir f:0 r:0.5
\move(54 39)
\fcir f:0 r:0.5
\move(57 0)
\fcir f:0 r:0.5
\move(57 3)
\fcir f:0 r:0.5
\move(57 6)
\fcir f:0 r:0.5
\move(57 9)
\fcir f:0 r:0.5
\move(57 12)
\fcir f:0 r:0.5
\move(57 15)
\fcir f:0 r:0.5
\move(57 18)
\fcir f:0 r:0.5
\move(57 21)
\fcir f:0 r:0.5
\move(57 24)
\fcir f:0 r:0.5
\move(57 27)
\fcir f:0 r:0.5
\move(57 30)
\fcir f:0 r:0.5
\move(57 33)
\fcir f:0 r:0.5
\move(57 36)
\fcir f:0 r:0.5
\move(57 39)
\fcir f:0 r:0.5
\move(60 0)
\fcir f:0 r:0.5
\move(60 3)
\fcir f:0 r:0.5
\move(60 6)
\fcir f:0 r:0.5
\move(60 9)
\fcir f:0 r:0.5
\move(60 12)
\fcir f:0 r:0.5
\move(60 15)
\fcir f:0 r:0.5
\move(60 18)
\fcir f:0 r:0.5
\move(60 21)
\fcir f:0 r:0.5
\move(60 24)
\fcir f:0 r:0.5
\move(60 27)
\fcir f:0 r:0.5
\move(60 30)
\fcir f:0 r:0.5
\move(60 33)
\fcir f:0 r:0.5
\move(60 36)
\fcir f:0 r:0.5
\move(60 39)
\fcir f:0 r:0.5
\move(63 0)
\fcir f:0 r:0.5
\move(63 3)
\fcir f:0 r:0.5
\move(63 6)
\fcir f:0 r:0.5
\move(63 9)
\fcir f:0 r:0.5
\move(63 12)
\fcir f:0 r:0.5
\move(63 15)
\fcir f:0 r:0.5
\move(63 18)
\fcir f:0 r:0.5
\move(63 21)
\fcir f:0 r:0.5
\move(63 24)
\fcir f:0 r:0.5
\move(63 27)
\fcir f:0 r:0.5
\move(63 30)
\fcir f:0 r:0.5
\move(63 33)
\fcir f:0 r:0.5
\move(63 36)
\fcir f:0 r:0.5
\move(63 39)
\fcir f:0 r:0.5
\move(66 0)
\fcir f:0 r:0.5
\move(66 3)
\fcir f:0 r:0.5
\move(66 6)
\fcir f:0 r:0.5
\move(66 9)
\fcir f:0 r:0.5
\move(66 12)
\fcir f:0 r:0.5
\move(66 15)
\fcir f:0 r:0.5
\move(66 18)
\fcir f:0 r:0.5
\move(66 21)
\fcir f:0 r:0.5
\move(66 24)
\fcir f:0 r:0.5
\move(66 27)
\fcir f:0 r:0.5
\move(66 30)
\fcir f:0 r:0.5
\move(66 33)
\fcir f:0 r:0.5
\move(66 36)
\fcir f:0 r:0.5
\move(66 39)
\fcir f:0 r:0.5
\move(69 0)
\fcir f:0 r:0.5
\move(69 3)
\fcir f:0 r:0.5
\move(69 6)
\fcir f:0 r:0.5
\move(69 9)
\fcir f:0 r:0.5
\move(69 12)
\fcir f:0 r:0.5
\move(69 15)
\fcir f:0 r:0.5
\move(69 18)
\fcir f:0 r:0.5
\move(69 21)
\fcir f:0 r:0.5
\move(69 24)
\fcir f:0 r:0.5
\move(69 27)
\fcir f:0 r:0.5
\move(69 30)
\fcir f:0 r:0.5
\move(69 33)
\fcir f:0 r:0.5
\move(69 36)
\fcir f:0 r:0.5
\move(69 39)
\fcir f:0 r:0.5
\move(72 0)
\fcir f:0 r:0.5
\move(72 3)
\fcir f:0 r:0.5
\move(72 6)
\fcir f:0 r:0.5
\move(72 9)
\fcir f:0 r:0.5
\move(72 12)
\fcir f:0 r:0.5
\move(72 15)
\fcir f:0 r:0.5
\move(72 18)
\fcir f:0 r:0.5
\move(72 21)
\fcir f:0 r:0.5
\move(72 24)
\fcir f:0 r:0.5
\move(72 27)
\fcir f:0 r:0.5
\move(72 30)
\fcir f:0 r:0.5
\move(72 33)
\fcir f:0 r:0.5
\move(72 36)
\fcir f:0 r:0.5
\move(72 39)
\fcir f:0 r:0.5
\move(75 0)
\fcir f:0 r:0.5
\move(75 3)
\fcir f:0 r:0.5
\move(75 6)
\fcir f:0 r:0.5
\move(75 9)
\fcir f:0 r:0.5
\move(75 12)
\fcir f:0 r:0.5
\move(75 15)
\fcir f:0 r:0.5
\move(75 18)
\fcir f:0 r:0.5
\move(75 21)
\fcir f:0 r:0.5
\move(75 24)
\fcir f:0 r:0.5
\move(75 27)
\fcir f:0 r:0.5
\move(75 30)
\fcir f:0 r:0.5
\move(75 33)
\fcir f:0 r:0.5
\move(75 36)
\fcir f:0 r:0.5
\move(75 39)
\fcir f:0 r:0.5
\move(78 0)
\fcir f:0 r:0.5
\move(78 3)
\fcir f:0 r:0.5
\move(78 6)
\fcir f:0 r:0.5
\move(78 9)
\fcir f:0 r:0.5
\move(78 12)
\fcir f:0 r:0.5
\move(78 15)
\fcir f:0 r:0.5
\move(78 18)
\fcir f:0 r:0.5
\move(78 21)
\fcir f:0 r:0.5
\move(78 24)
\fcir f:0 r:0.5
\move(78 27)
\fcir f:0 r:0.5
\move(78 30)
\fcir f:0 r:0.5
\move(78 33)
\fcir f:0 r:0.5
\move(78 36)
\fcir f:0 r:0.5
\move(78 39)
\fcir f:0 r:0.5
\move(81 0)
\fcir f:0 r:0.5
\move(81 3)
\fcir f:0 r:0.5
\move(81 6)
\fcir f:0 r:0.5
\move(81 9)
\fcir f:0 r:0.5
\move(81 12)
\fcir f:0 r:0.5
\move(81 15)
\fcir f:0 r:0.5
\move(81 18)
\fcir f:0 r:0.5
\move(81 21)
\fcir f:0 r:0.5
\move(81 24)
\fcir f:0 r:0.5
\move(81 27)
\fcir f:0 r:0.5
\move(81 30)
\fcir f:0 r:0.5
\move(81 33)
\fcir f:0 r:0.5
\move(81 36)
\fcir f:0 r:0.5
\move(81 39)
\fcir f:0 r:0.5
\move(84 0)
\fcir f:0 r:0.5
\move(84 3)
\fcir f:0 r:0.5
\move(84 6)
\fcir f:0 r:0.5
\move(84 9)
\fcir f:0 r:0.5
\move(84 12)
\fcir f:0 r:0.5
\move(84 15)
\fcir f:0 r:0.5
\move(84 18)
\fcir f:0 r:0.5
\move(84 21)
\fcir f:0 r:0.5
\move(84 24)
\fcir f:0 r:0.5
\move(84 27)
\fcir f:0 r:0.5
\move(84 30)
\fcir f:0 r:0.5
\move(84 33)
\fcir f:0 r:0.5
\move(84 36)
\fcir f:0 r:0.5
\move(84 39)
\fcir f:0 r:0.5
\move(87 0)
\fcir f:0 r:0.5
\move(87 3)
\fcir f:0 r:0.5
\move(87 6)
\fcir f:0 r:0.5
\move(87 9)
\fcir f:0 r:0.5
\move(87 12)
\fcir f:0 r:0.5
\move(87 15)
\fcir f:0 r:0.5
\move(87 18)
\fcir f:0 r:0.5
\move(87 21)
\fcir f:0 r:0.5
\move(87 24)
\fcir f:0 r:0.5
\move(87 27)
\fcir f:0 r:0.5
\move(87 30)
\fcir f:0 r:0.5
\move(87 33)
\fcir f:0 r:0.5
\move(87 36)
\fcir f:0 r:0.5
\move(87 39)
\fcir f:0 r:0.5
\move(90 0)
\fcir f:0 r:0.5
\move(90 3)
\fcir f:0 r:0.5
\move(90 6)
\fcir f:0 r:0.5
\move(90 9)
\fcir f:0 r:0.5
\move(90 12)
\fcir f:0 r:0.5
\move(90 15)
\fcir f:0 r:0.5
\move(90 18)
\fcir f:0 r:0.5
\move(90 21)
\fcir f:0 r:0.5
\move(90 24)
\fcir f:0 r:0.5
\move(90 27)
\fcir f:0 r:0.5
\move(90 30)
\fcir f:0 r:0.5
\move(90 33)
\fcir f:0 r:0.5
\move(90 36)
\fcir f:0 r:0.5
\move(90 39)
\fcir f:0 r:0.5
\move(93 0)
\fcir f:0 r:0.5
\move(93 3)
\fcir f:0 r:0.5
\move(93 6)
\fcir f:0 r:0.5
\move(93 9)
\fcir f:0 r:0.5
\move(93 12)
\fcir f:0 r:0.5
\move(93 15)
\fcir f:0 r:0.5
\move(93 18)
\fcir f:0 r:0.5
\move(93 21)
\fcir f:0 r:0.5
\move(93 24)
\fcir f:0 r:0.5
\move(93 27)
\fcir f:0 r:0.5
\move(93 30)
\fcir f:0 r:0.5
\move(93 33)
\fcir f:0 r:0.5
\move(93 36)
\fcir f:0 r:0.5
\move(93 39)
\fcir f:0 r:0.5
\move(96 0)
\fcir f:0 r:0.5
\move(96 3)
\fcir f:0 r:0.5
\move(96 6)
\fcir f:0 r:0.5
\move(96 9)
\fcir f:0 r:0.5
\move(96 12)
\fcir f:0 r:0.5
\move(96 15)
\fcir f:0 r:0.5
\move(96 18)
\fcir f:0 r:0.5
\move(96 21)
\fcir f:0 r:0.5
\move(96 24)
\fcir f:0 r:0.5
\move(96 27)
\fcir f:0 r:0.5
\move(96 30)
\fcir f:0 r:0.5
\move(96 33)
\fcir f:0 r:0.5
\move(96 36)
\fcir f:0 r:0.5
\move(96 39)
\fcir f:0 r:0.5
\move(99 0)
\fcir f:0 r:0.5
\move(99 3)
\fcir f:0 r:0.5
\move(99 6)
\fcir f:0 r:0.5
\move(99 9)
\fcir f:0 r:0.5
\move(99 12)
\fcir f:0 r:0.5
\move(99 15)
\fcir f:0 r:0.5
\move(99 18)
\fcir f:0 r:0.5
\move(99 21)
\fcir f:0 r:0.5
\move(99 24)
\fcir f:0 r:0.5
\move(99 27)
\fcir f:0 r:0.5
\move(99 30)
\fcir f:0 r:0.5
\move(99 33)
\fcir f:0 r:0.5
\move(99 36)
\fcir f:0 r:0.5
\move(99 39)
\fcir f:0 r:0.5
\move(102 0)
\fcir f:0 r:0.5
\move(102 3)
\fcir f:0 r:0.5
\move(102 6)
\fcir f:0 r:0.5
\move(102 9)
\fcir f:0 r:0.5
\move(102 12)
\fcir f:0 r:0.5
\move(102 15)
\fcir f:0 r:0.5
\move(102 18)
\fcir f:0 r:0.5
\move(102 21)
\fcir f:0 r:0.5
\move(102 24)
\fcir f:0 r:0.5
\move(102 27)
\fcir f:0 r:0.5
\move(102 30)
\fcir f:0 r:0.5
\move(102 33)
\fcir f:0 r:0.5
\move(102 36)
\fcir f:0 r:0.5
\move(102 39)
\fcir f:0 r:0.5
\move(105 0)
\fcir f:0 r:0.5
\move(105 3)
\fcir f:0 r:0.5
\move(105 6)
\fcir f:0 r:0.5
\move(105 9)
\fcir f:0 r:0.5
\move(105 12)
\fcir f:0 r:0.5
\move(105 15)
\fcir f:0 r:0.5
\move(105 18)
\fcir f:0 r:0.5
\move(105 21)
\fcir f:0 r:0.5
\move(105 24)
\fcir f:0 r:0.5
\move(105 27)
\fcir f:0 r:0.5
\move(105 30)
\fcir f:0 r:0.5
\move(105 33)
\fcir f:0 r:0.5
\move(105 36)
\fcir f:0 r:0.5
\move(105 39)
\fcir f:0 r:0.5
\move(108 0)
\fcir f:0 r:0.5
\move(108 3)
\fcir f:0 r:0.5
\move(108 6)
\fcir f:0 r:0.5
\move(108 9)
\fcir f:0 r:0.5
\move(108 12)
\fcir f:0 r:0.5
\move(108 15)
\fcir f:0 r:0.5
\move(108 18)
\fcir f:0 r:0.5
\move(108 21)
\fcir f:0 r:0.5
\move(108 24)
\fcir f:0 r:0.5
\move(108 27)
\fcir f:0 r:0.5
\move(108 30)
\fcir f:0 r:0.5
\move(108 33)
\fcir f:0 r:0.5
\move(108 36)
\fcir f:0 r:0.5
\move(108 39)
\fcir f:0 r:0.5
\move(111 0)
\fcir f:0 r:0.5
\move(111 3)
\fcir f:0 r:0.5
\move(111 6)
\fcir f:0 r:0.5
\move(111 9)
\fcir f:0 r:0.5
\move(111 12)
\fcir f:0 r:0.5
\move(111 15)
\fcir f:0 r:0.5
\move(111 18)
\fcir f:0 r:0.5
\move(111 21)
\fcir f:0 r:0.5
\move(111 24)
\fcir f:0 r:0.5
\move(111 27)
\fcir f:0 r:0.5
\move(111 30)
\fcir f:0 r:0.5
\move(111 33)
\fcir f:0 r:0.5
\move(111 36)
\fcir f:0 r:0.5
\move(111 39)
\fcir f:0 r:0.5
\move(114 0)
\fcir f:0 r:0.5
\move(114 3)
\fcir f:0 r:0.5
\move(114 6)
\fcir f:0 r:0.5
\move(114 9)
\fcir f:0 r:0.5
\move(114 12)
\fcir f:0 r:0.5
\move(114 15)
\fcir f:0 r:0.5
\move(114 18)
\fcir f:0 r:0.5
\move(114 21)
\fcir f:0 r:0.5
\move(114 24)
\fcir f:0 r:0.5
\move(114 27)
\fcir f:0 r:0.5
\move(114 30)
\fcir f:0 r:0.5
\move(114 33)
\fcir f:0 r:0.5
\move(114 36)
\fcir f:0 r:0.5
\move(114 39)
\fcir f:0 r:0.5
\move(117 0)
\fcir f:0 r:0.5
\move(117 3)
\fcir f:0 r:0.5
\move(117 6)
\fcir f:0 r:0.5
\move(117 9)
\fcir f:0 r:0.5
\move(117 12)
\fcir f:0 r:0.5
\move(117 15)
\fcir f:0 r:0.5
\move(117 18)
\fcir f:0 r:0.5
\move(117 21)
\fcir f:0 r:0.5
\move(117 24)
\fcir f:0 r:0.5
\move(117 27)
\fcir f:0 r:0.5
\move(117 30)
\fcir f:0 r:0.5
\move(117 33)
\fcir f:0 r:0.5
\move(117 36)
\fcir f:0 r:0.5
\move(117 39)
\fcir f:0 r:0.5
\move(120 0)
\fcir f:0 r:0.5
\move(120 3)
\fcir f:0 r:0.5
\move(120 6)
\fcir f:0 r:0.5
\move(120 9)
\fcir f:0 r:0.5
\move(120 12)
\fcir f:0 r:0.5
\move(120 15)
\fcir f:0 r:0.5
\move(120 18)
\fcir f:0 r:0.5
\move(120 21)
\fcir f:0 r:0.5
\move(120 24)
\fcir f:0 r:0.5
\move(120 27)
\fcir f:0 r:0.5
\move(120 30)
\fcir f:0 r:0.5
\move(120 33)
\fcir f:0 r:0.5
\move(120 36)
\fcir f:0 r:0.5
\move(120 39)
\fcir f:0 r:0.5
\move(123 0)
\fcir f:0 r:0.5
\move(123 3)
\fcir f:0 r:0.5
\move(123 6)
\fcir f:0 r:0.5
\move(123 9)
\fcir f:0 r:0.5
\move(123 12)
\fcir f:0 r:0.5
\move(123 15)
\fcir f:0 r:0.5
\move(123 18)
\fcir f:0 r:0.5
\move(123 21)
\fcir f:0 r:0.5
\move(123 24)
\fcir f:0 r:0.5
\move(123 27)
\fcir f:0 r:0.5
\move(123 30)
\fcir f:0 r:0.5
\move(123 33)
\fcir f:0 r:0.5
\move(123 36)
\fcir f:0 r:0.5
\move(123 39)
\fcir f:0 r:0.5
\move(126 0)
\fcir f:0 r:0.5
\move(126 3)
\fcir f:0 r:0.5
\move(126 6)
\fcir f:0 r:0.5
\move(126 9)
\fcir f:0 r:0.5
\move(126 12)
\fcir f:0 r:0.5
\move(126 15)
\fcir f:0 r:0.5
\move(126 18)
\fcir f:0 r:0.5
\move(126 21)
\fcir f:0 r:0.5
\move(126 24)
\fcir f:0 r:0.5
\move(126 27)
\fcir f:0 r:0.5
\move(126 30)
\fcir f:0 r:0.5
\move(126 33)
\fcir f:0 r:0.5
\move(126 36)
\fcir f:0 r:0.5
\move(126 39)
\fcir f:0 r:0.5
\move(129 0)
\fcir f:0 r:0.5
\move(129 3)
\fcir f:0 r:0.5
\move(129 6)
\fcir f:0 r:0.5
\move(129 9)
\fcir f:0 r:0.5
\move(129 12)
\fcir f:0 r:0.5
\move(129 15)
\fcir f:0 r:0.5
\move(129 18)
\fcir f:0 r:0.5
\move(129 21)
\fcir f:0 r:0.5
\move(129 24)
\fcir f:0 r:0.5
\move(129 27)
\fcir f:0 r:0.5
\move(129 30)
\fcir f:0 r:0.5
\move(129 33)
\fcir f:0 r:0.5
\move(129 36)
\fcir f:0 r:0.5
\move(129 39)
\fcir f:0 r:0.5
\move(132 0)
\fcir f:0 r:0.5
\move(132 3)
\fcir f:0 r:0.5
\move(132 6)
\fcir f:0 r:0.5
\move(132 9)
\fcir f:0 r:0.5
\move(132 12)
\fcir f:0 r:0.5
\move(132 15)
\fcir f:0 r:0.5
\move(132 18)
\fcir f:0 r:0.5
\move(132 21)
\fcir f:0 r:0.5
\move(132 24)
\fcir f:0 r:0.5
\move(132 27)
\fcir f:0 r:0.5
\move(132 30)
\fcir f:0 r:0.5
\move(132 33)
\fcir f:0 r:0.5
\move(132 36)
\fcir f:0 r:0.5
\move(132 39)
\fcir f:0 r:0.5
\move(135 0)
\fcir f:0 r:0.5
\move(135 3)
\fcir f:0 r:0.5
\move(135 6)
\fcir f:0 r:0.5
\move(135 9)
\fcir f:0 r:0.5
\move(135 12)
\fcir f:0 r:0.5
\move(135 15)
\fcir f:0 r:0.5
\move(135 18)
\fcir f:0 r:0.5
\move(135 21)
\fcir f:0 r:0.5
\move(135 24)
\fcir f:0 r:0.5
\move(135 27)
\fcir f:0 r:0.5
\move(135 30)
\fcir f:0 r:0.5
\move(135 33)
\fcir f:0 r:0.5
\move(135 36)
\fcir f:0 r:0.5
\move(135 39)
\fcir f:0 r:0.5
\move(138 0)
\fcir f:0 r:0.5
\move(138 3)
\fcir f:0 r:0.5
\move(138 6)
\fcir f:0 r:0.5
\move(138 9)
\fcir f:0 r:0.5
\move(138 12)
\fcir f:0 r:0.5
\move(138 15)
\fcir f:0 r:0.5
\move(138 18)
\fcir f:0 r:0.5
\move(138 21)
\fcir f:0 r:0.5
\move(138 24)
\fcir f:0 r:0.5
\move(138 27)
\fcir f:0 r:0.5
\move(138 30)
\fcir f:0 r:0.5
\move(138 33)
\fcir f:0 r:0.5
\move(138 36)
\fcir f:0 r:0.5
\move(138 39)
\fcir f:0 r:0.5
\move(141 0)
\fcir f:0 r:0.5
\move(141 3)
\fcir f:0 r:0.5
\move(141 6)
\fcir f:0 r:0.5
\move(141 9)
\fcir f:0 r:0.5
\move(141 12)
\fcir f:0 r:0.5
\move(141 15)
\fcir f:0 r:0.5
\move(141 18)
\fcir f:0 r:0.5
\move(141 21)
\fcir f:0 r:0.5
\move(141 24)
\fcir f:0 r:0.5
\move(141 27)
\fcir f:0 r:0.5
\move(141 30)
\fcir f:0 r:0.5
\move(141 33)
\fcir f:0 r:0.5
\move(141 36)
\fcir f:0 r:0.5
\move(141 39)
\fcir f:0 r:0.5
\move(144 0)
\fcir f:0 r:0.5
\move(144 3)
\fcir f:0 r:0.5
\move(144 6)
\fcir f:0 r:0.5
\move(144 9)
\fcir f:0 r:0.5
\move(144 12)
\fcir f:0 r:0.5
\move(144 15)
\fcir f:0 r:0.5
\move(144 18)
\fcir f:0 r:0.5
\move(144 21)
\fcir f:0 r:0.5
\move(144 24)
\fcir f:0 r:0.5
\move(144 27)
\fcir f:0 r:0.5
\move(144 30)
\fcir f:0 r:0.5
\move(144 33)
\fcir f:0 r:0.5
\move(144 36)
\fcir f:0 r:0.5
\move(144 39)
\fcir f:0 r:0.5
\move(147 0)
\fcir f:0 r:0.5
\move(147 3)
\fcir f:0 r:0.5
\move(147 6)
\fcir f:0 r:0.5
\move(147 9)
\fcir f:0 r:0.5
\move(147 12)
\fcir f:0 r:0.5
\move(147 15)
\fcir f:0 r:0.5
\move(147 18)
\fcir f:0 r:0.5
\move(147 21)
\fcir f:0 r:0.5
\move(147 24)
\fcir f:0 r:0.5
\move(147 27)
\fcir f:0 r:0.5
\move(147 30)
\fcir f:0 r:0.5
\move(147 33)
\fcir f:0 r:0.5
\move(147 36)
\fcir f:0 r:0.5
\move(147 39)
\fcir f:0 r:0.5
\move(150 0)
\fcir f:0 r:0.5
\move(150 3)
\fcir f:0 r:0.5
\move(150 6)
\fcir f:0 r:0.5
\move(150 9)
\fcir f:0 r:0.5
\move(150 12)
\fcir f:0 r:0.5
\move(150 15)
\fcir f:0 r:0.5
\move(150 18)
\fcir f:0 r:0.5
\move(150 21)
\fcir f:0 r:0.5
\move(150 24)
\fcir f:0 r:0.5
\move(150 27)
\fcir f:0 r:0.5
\move(150 30)
\fcir f:0 r:0.5
\move(150 33)
\fcir f:0 r:0.5
\move(150 36)
\fcir f:0 r:0.5
\move(150 39)
\fcir f:0 r:0.5
\move(153 0)
\fcir f:0 r:0.5
\move(153 3)
\fcir f:0 r:0.5
\move(153 6)
\fcir f:0 r:0.5
\move(153 9)
\fcir f:0 r:0.5
\move(153 12)
\fcir f:0 r:0.5
\move(153 15)
\fcir f:0 r:0.5
\move(153 18)
\fcir f:0 r:0.5
\move(153 21)
\fcir f:0 r:0.5
\move(153 24)
\fcir f:0 r:0.5
\move(153 27)
\fcir f:0 r:0.5
\move(153 30)
\fcir f:0 r:0.5
\move(153 33)
\fcir f:0 r:0.5
\move(153 36)
\fcir f:0 r:0.5
\move(153 39)
\fcir f:0 r:0.5
\move(156 0)
\fcir f:0 r:0.5
\move(156 3)
\fcir f:0 r:0.5
\move(156 6)
\fcir f:0 r:0.5
\move(156 9)
\fcir f:0 r:0.5
\move(156 12)
\fcir f:0 r:0.5
\move(156 15)
\fcir f:0 r:0.5
\move(156 18)
\fcir f:0 r:0.5
\move(156 21)
\fcir f:0 r:0.5
\move(156 24)
\fcir f:0 r:0.5
\move(156 27)
\fcir f:0 r:0.5
\move(156 30)
\fcir f:0 r:0.5
\move(156 33)
\fcir f:0 r:0.5
\move(156 36)
\fcir f:0 r:0.5
\move(156 39)
\fcir f:0 r:0.5
\move(159 0)
\fcir f:0 r:0.5
\move(159 3)
\fcir f:0 r:0.5
\move(159 6)
\fcir f:0 r:0.5
\move(159 9)
\fcir f:0 r:0.5
\move(159 12)
\fcir f:0 r:0.5
\move(159 15)
\fcir f:0 r:0.5
\move(159 18)
\fcir f:0 r:0.5
\move(159 21)
\fcir f:0 r:0.5
\move(159 24)
\fcir f:0 r:0.5
\move(159 27)
\fcir f:0 r:0.5
\move(159 30)
\fcir f:0 r:0.5
\move(159 33)
\fcir f:0 r:0.5
\move(159 36)
\fcir f:0 r:0.5
\move(159 39)
\fcir f:0 r:0.5
\move(162 0)
\fcir f:0 r:0.5
\move(162 3)
\fcir f:0 r:0.5
\move(162 6)
\fcir f:0 r:0.5
\move(162 9)
\fcir f:0 r:0.5
\move(162 12)
\fcir f:0 r:0.5
\move(162 15)
\fcir f:0 r:0.5
\move(162 18)
\fcir f:0 r:0.5
\move(162 21)
\fcir f:0 r:0.5
\move(162 24)
\fcir f:0 r:0.5
\move(162 27)
\fcir f:0 r:0.5
\move(162 30)
\fcir f:0 r:0.5
\move(162 33)
\fcir f:0 r:0.5
\move(162 36)
\fcir f:0 r:0.5
\move(162 39)
\fcir f:0 r:0.5
\move(165 0)
\fcir f:0 r:0.5
\move(165 3)
\fcir f:0 r:0.5
\move(165 6)
\fcir f:0 r:0.5
\move(165 9)
\fcir f:0 r:0.5
\move(165 12)
\fcir f:0 r:0.5
\move(165 15)
\fcir f:0 r:0.5
\move(165 18)
\fcir f:0 r:0.5
\move(165 21)
\fcir f:0 r:0.5
\move(165 24)
\fcir f:0 r:0.5
\move(165 27)
\fcir f:0 r:0.5
\move(165 30)
\fcir f:0 r:0.5
\move(165 33)
\fcir f:0 r:0.5
\move(165 36)
\fcir f:0 r:0.5
\move(165 39)
\fcir f:0 r:0.5
\move(168 0)
\fcir f:0 r:0.5
\move(168 3)
\fcir f:0 r:0.5
\move(168 6)
\fcir f:0 r:0.5
\move(168 9)
\fcir f:0 r:0.5
\move(168 12)
\fcir f:0 r:0.5
\move(168 15)
\fcir f:0 r:0.5
\move(168 18)
\fcir f:0 r:0.5
\move(168 21)
\fcir f:0 r:0.5
\move(168 24)
\fcir f:0 r:0.5
\move(168 27)
\fcir f:0 r:0.5
\move(168 30)
\fcir f:0 r:0.5
\move(168 33)
\fcir f:0 r:0.5
\move(168 36)
\fcir f:0 r:0.5
\move(168 39)
\fcir f:0 r:0.5
\move(171 0)
\fcir f:0 r:0.5
\move(171 3)
\fcir f:0 r:0.5
\move(171 6)
\fcir f:0 r:0.5
\move(171 9)
\fcir f:0 r:0.5
\move(171 12)
\fcir f:0 r:0.5
\move(171 15)
\fcir f:0 r:0.5
\move(171 18)
\fcir f:0 r:0.5
\move(171 21)
\fcir f:0 r:0.5
\move(171 24)
\fcir f:0 r:0.5
\move(171 27)
\fcir f:0 r:0.5
\move(171 30)
\fcir f:0 r:0.5
\move(171 33)
\fcir f:0 r:0.5
\move(171 36)
\fcir f:0 r:0.5
\move(171 39)
\fcir f:0 r:0.5
\arrowheadtype t:V
\move(0 0)
\avec(180 0)
\move(0 0)
\avec(0 48)
\htext(182 0){$\mathbb{N}$}
\htext(-13 42){$\mathbb{N}$}
\move(0 0)
\fcir f:0 r:1.5
\move(3 0)
\fcir f:0 r:1.5
\move(6 0)
\fcir f:0 r:1.5
\move(9 0)
\fcir f:0 r:1.5
\move(12 0)
\fcir f:0 r:1.5
\move(12 3)
\fcir f:0 r:1.5
\move(15 0)
\fcir f:0 r:1.5
\move(15 3)
\fcir f:0 r:1.5
\move(18 0)
\fcir f:0 r:1.5
\move(18 3)
\fcir f:0 r:1.5
\move(21 0)
\fcir f:0 r:1.5
\move(21 3)
\fcir f:0 r:1.5
\move(24 0)
\fcir f:0 r:1.5
\move(24 3)
\fcir f:0 r:1.5
\move(24 6)
\fcir f:0 r:1.5
\move(27 0)
\fcir f:0 r:1.5
\move(27 3)
\fcir f:0 r:1.5
\move(27 6)
\fcir f:0 r:1.5
\move(30 0)
\fcir f:0 r:1.5
\move(30 3)
\fcir f:0 r:1.5
\move(30 6)
\fcir f:0 r:1.5
\move(33 0)
\fcir f:0 r:1.5
\move(33 3)
\fcir f:0 r:1.5
\move(33 6)
\fcir f:0 r:1.5
\move(36 0)
\fcir f:0 r:1.5
\move(36 3)
\fcir f:0 r:1.5
\move(36 6)
\fcir f:0 r:1.5
\move(36 9)
\fcir f:0 r:1.5
\move(39 0)
\fcir f:0 r:1.5
\move(39 3)
\fcir f:0 r:1.5
\move(39 6)
\fcir f:0 r:1.5
\move(39 9)
\fcir f:0 r:1.5
\move(42 0)
\fcir f:0 r:1.5
\move(42 3)
\fcir f:0 r:1.5
\move(42 6)
\fcir f:0 r:1.5
\move(42 9)
\fcir f:0 r:1.5
\move(45 0)
\fcir f:0 r:1.5
\move(45 3)
\fcir f:0 r:1.5
\move(45 6)
\fcir f:0 r:1.5
\move(45 9)
\fcir f:0 r:1.5
\move(48 0)
\fcir f:0 r:1.5
\move(48 3)
\fcir f:0 r:1.5
\move(48 6)
\fcir f:0 r:1.5
\move(48 9)
\fcir f:0 r:1.5
\move(48 12)
\fcir f:0 r:1.5
\move(51 0)
\fcir f:0 r:1.5
\move(51 3)
\fcir f:0 r:1.5
\move(51 6)
\fcir f:0 r:1.5
\move(51 9)
\fcir f:0 r:1.5
\move(51 12)
\fcir f:0 r:1.5
\move(54 0)
\fcir f:0 r:1.5
\move(54 3)
\fcir f:0 r:1.5
\move(54 6)
\fcir f:0 r:1.5
\move(54 9)
\fcir f:0 r:1.5
\move(54 12)
\fcir f:0 r:1.5
\move(57 0)
\fcir f:0 r:1.5
\move(57 3)
\fcir f:0 r:1.5
\move(57 6)
\fcir f:0 r:1.5
\move(57 9)
\fcir f:0 r:1.5
\move(57 12)
\fcir f:0 r:1.5
\move(60 0)
\fcir f:0 r:1.5
\move(60 3)
\fcir f:0 r:1.5
\move(60 6)
\fcir f:0 r:1.5
\move(60 9)
\fcir f:0 r:1.5
\move(60 12)
\fcir f:0 r:1.5
\move(60 15)
\fcir f:0 r:1.5
\move(63 0)
\fcir f:0 r:1.5
\move(63 3)
\fcir f:0 r:1.5
\move(63 6)
\fcir f:0 r:1.5
\move(63 9)
\fcir f:0 r:1.5
\move(63 12)
\fcir f:0 r:1.5
\move(63 15)
\fcir f:0 r:1.5
\move(66 0)
\fcir f:0 r:1.5
\move(66 3)
\fcir f:0 r:1.5
\move(66 6)
\fcir f:0 r:1.5
\move(66 9)
\fcir f:0 r:1.5
\move(66 12)
\fcir f:0 r:1.5
\move(66 15)
\fcir f:0 r:1.5
\move(69 0)
\fcir f:0 r:1.5
\move(69 3)
\fcir f:0 r:1.5
\move(69 6)
\fcir f:0 r:1.5
\move(69 9)
\fcir f:0 r:1.5
\move(69 12)
\fcir f:0 r:1.5
\move(69 15)
\fcir f:0 r:1.5
\move(72 0)
\fcir f:0 r:1.5
\move(72 3)
\fcir f:0 r:1.5
\move(72 6)
\fcir f:0 r:1.5
\move(72 9)
\fcir f:0 r:1.5
\move(72 12)
\fcir f:0 r:1.5
\move(72 15)
\fcir f:0 r:1.5
\move(72 18)
\fcir f:0 r:1.5
\move(75 0)
\fcir f:0 r:1.5
\move(75 3)
\fcir f:0 r:1.5
\move(75 6)
\fcir f:0 r:1.5
\move(75 9)
\fcir f:0 r:1.5
\move(75 12)
\fcir f:0 r:1.5
\move(75 15)
\fcir f:0 r:1.5
\move(75 18)
\fcir f:0 r:1.5
\move(78 0)
\fcir f:0 r:1.5
\move(78 3)
\fcir f:0 r:1.5
\move(78 6)
\fcir f:0 r:1.5
\move(78 9)
\fcir f:0 r:1.5
\move(78 12)
\fcir f:0 r:1.5
\move(78 15)
\fcir f:0 r:1.5
\move(78 18)
\fcir f:0 r:1.5
\move(81 0)
\fcir f:0 r:1.5
\move(81 3)
\fcir f:0 r:1.5
\move(81 6)
\fcir f:0 r:1.5
\move(81 9)
\fcir f:0 r:1.5
\move(81 12)
\fcir f:0 r:1.5
\move(81 15)
\fcir f:0 r:1.5
\move(81 18)
\fcir f:0 r:1.5
\move(84 0)
\fcir f:0 r:1.5
\move(84 3)
\fcir f:0 r:1.5
\move(84 6)
\fcir f:0 r:1.5
\move(84 9)
\fcir f:0 r:1.5
\move(84 12)
\fcir f:0 r:1.5
\move(84 15)
\fcir f:0 r:1.5
\move(84 18)
\fcir f:0 r:1.5
\move(84 21)
\fcir f:0 r:1.5
\move(87 0)
\fcir f:0 r:1.5
\move(87 3)
\fcir f:0 r:1.5
\move(87 6)
\fcir f:0 r:1.5
\move(87 9)
\fcir f:0 r:1.5
\move(87 12)
\fcir f:0 r:1.5
\move(87 15)
\fcir f:0 r:1.5
\move(87 18)
\fcir f:0 r:1.5
\move(87 21)
\fcir f:0 r:1.5
\move(90 0)
\fcir f:0 r:1.5
\move(90 3)
\fcir f:0 r:1.5
\move(90 6)
\fcir f:0 r:1.5
\move(90 9)
\fcir f:0 r:1.5
\move(90 12)
\fcir f:0 r:1.5
\move(90 15)
\fcir f:0 r:1.5
\move(90 18)
\fcir f:0 r:1.5
\move(90 21)
\fcir f:0 r:1.5
\move(93 0)
\fcir f:0 r:1.5
\move(93 3)
\fcir f:0 r:1.5
\move(93 6)
\fcir f:0 r:1.5
\move(93 9)
\fcir f:0 r:1.5
\move(93 12)
\fcir f:0 r:1.5
\move(93 15)
\fcir f:0 r:1.5
\move(93 18)
\fcir f:0 r:1.5
\move(93 21)
\fcir f:0 r:1.5
\move(96 0)
\fcir f:0 r:1.5
\move(96 3)
\fcir f:0 r:1.5
\move(96 6)
\fcir f:0 r:1.5
\move(96 9)
\fcir f:0 r:1.5
\move(96 12)
\fcir f:0 r:1.5
\move(96 15)
\fcir f:0 r:1.5
\move(96 18)
\fcir f:0 r:1.5
\move(96 21)
\fcir f:0 r:1.5
\move(96 24)
\fcir f:0 r:1.5
\move(99 0)
\fcir f:0 r:1.5
\move(99 3)
\fcir f:0 r:1.5
\move(99 6)
\fcir f:0 r:1.5
\move(99 9)
\fcir f:0 r:1.5
\move(99 12)
\fcir f:0 r:1.5
\move(99 15)
\fcir f:0 r:1.5
\move(99 18)
\fcir f:0 r:1.5
\move(99 21)
\fcir f:0 r:1.5
\move(99 24)
\fcir f:0 r:1.5
\move(102 0)
\fcir f:0 r:1.5
\move(102 3)
\fcir f:0 r:1.5
\move(102 6)
\fcir f:0 r:1.5
\move(102 9)
\fcir f:0 r:1.5
\move(102 12)
\fcir f:0 r:1.5
\move(102 15)
\fcir f:0 r:1.5
\move(102 18)
\fcir f:0 r:1.5
\move(102 21)
\fcir f:0 r:1.5
\move(102 24)
\fcir f:0 r:1.5
\move(105 0)
\fcir f:0 r:1.5
\move(105 3)
\fcir f:0 r:1.5
\move(105 6)
\fcir f:0 r:1.5
\move(105 9)
\fcir f:0 r:1.5
\move(105 12)
\fcir f:0 r:1.5
\move(105 15)
\fcir f:0 r:1.5
\move(105 18)
\fcir f:0 r:1.5
\move(105 21)
\fcir f:0 r:1.5
\move(105 24)
\fcir f:0 r:1.5
\move(108 0)
\fcir f:0 r:1.5
\move(108 3)
\fcir f:0 r:1.5
\move(108 6)
\fcir f:0 r:1.5
\move(108 9)
\fcir f:0 r:1.5
\move(108 12)
\fcir f:0 r:1.5
\move(108 15)
\fcir f:0 r:1.5
\move(108 18)
\fcir f:0 r:1.5
\move(108 21)
\fcir f:0 r:1.5
\move(108 24)
\fcir f:0 r:1.5
\move(108 27)
\fcir f:0 r:1.5
\move(111 0)
\fcir f:0 r:1.5
\move(111 3)
\fcir f:0 r:1.5
\move(111 6)
\fcir f:0 r:1.5
\move(111 9)
\fcir f:0 r:1.5
\move(111 12)
\fcir f:0 r:1.5
\move(111 15)
\fcir f:0 r:1.5
\move(111 18)
\fcir f:0 r:1.5
\move(111 21)
\fcir f:0 r:1.5
\move(111 24)
\fcir f:0 r:1.5
\move(111 27)
\fcir f:0 r:1.5
\move(114 0)
\fcir f:0 r:1.5
\move(114 3)
\fcir f:0 r:1.5
\move(114 6)
\fcir f:0 r:1.5
\move(114 9)
\fcir f:0 r:1.5
\move(114 12)
\fcir f:0 r:1.5
\move(114 15)
\fcir f:0 r:1.5
\move(114 18)
\fcir f:0 r:1.5
\move(114 21)
\fcir f:0 r:1.5
\move(114 24)
\fcir f:0 r:1.5
\move(114 27)
\fcir f:0 r:1.5
\move(117 0)
\fcir f:0 r:1.5
\move(117 3)
\fcir f:0 r:1.5
\move(117 6)
\fcir f:0 r:1.5
\move(117 9)
\fcir f:0 r:1.5
\move(117 12)
\fcir f:0 r:1.5
\move(117 15)
\fcir f:0 r:1.5
\move(117 18)
\fcir f:0 r:1.5
\move(117 21)
\fcir f:0 r:1.5
\move(117 24)
\fcir f:0 r:1.5
\move(117 27)
\fcir f:0 r:1.5
\move(120 0)
\fcir f:0 r:1.5
\move(120 3)
\fcir f:0 r:1.5
\move(120 6)
\fcir f:0 r:1.5
\move(120 9)
\fcir f:0 r:1.5
\move(120 12)
\fcir f:0 r:1.5
\move(120 15)
\fcir f:0 r:1.5
\move(120 18)
\fcir f:0 r:1.5
\move(120 21)
\fcir f:0 r:1.5
\move(120 24)
\fcir f:0 r:1.5
\move(120 27)
\fcir f:0 r:1.5
\move(120 30)
\fcir f:0 r:1.5
\move(123 0)
\fcir f:0 r:1.5
\move(123 3)
\fcir f:0 r:1.5
\move(123 6)
\fcir f:0 r:1.5
\move(123 9)
\fcir f:0 r:1.5
\move(123 12)
\fcir f:0 r:1.5
\move(123 15)
\fcir f:0 r:1.5
\move(123 18)
\fcir f:0 r:1.5
\move(123 21)
\fcir f:0 r:1.5
\move(123 24)
\fcir f:0 r:1.5
\move(123 27)
\fcir f:0 r:1.5
\move(123 30)
\fcir f:0 r:1.5
\move(126 0)
\fcir f:0 r:1.5
\move(126 3)
\fcir f:0 r:1.5
\move(126 6)
\fcir f:0 r:1.5
\move(126 9)
\fcir f:0 r:1.5
\move(126 12)
\fcir f:0 r:1.5
\move(126 15)
\fcir f:0 r:1.5
\move(126 18)
\fcir f:0 r:1.5
\move(126 21)
\fcir f:0 r:1.5
\move(126 24)
\fcir f:0 r:1.5
\move(126 27)
\fcir f:0 r:1.5
\move(126 30)
\fcir f:0 r:1.5
\move(129 0)
\fcir f:0 r:1.5
\move(129 3)
\fcir f:0 r:1.5
\move(129 6)
\fcir f:0 r:1.5
\move(129 9)
\fcir f:0 r:1.5
\move(129 12)
\fcir f:0 r:1.5
\move(129 15)
\fcir f:0 r:1.5
\move(129 18)
\fcir f:0 r:1.5
\move(129 21)
\fcir f:0 r:1.5
\move(129 24)
\fcir f:0 r:1.5
\move(129 27)
\fcir f:0 r:1.5
\move(129 30)
\fcir f:0 r:1.5
\move(132 0)
\fcir f:0 r:1.5
\move(132 3)
\fcir f:0 r:1.5
\move(132 6)
\fcir f:0 r:1.5
\move(132 9)
\fcir f:0 r:1.5
\move(132 12)
\fcir f:0 r:1.5
\move(132 15)
\fcir f:0 r:1.5
\move(132 18)
\fcir f:0 r:1.5
\move(132 21)
\fcir f:0 r:1.5
\move(132 24)
\fcir f:0 r:1.5
\move(132 27)
\fcir f:0 r:1.5
\move(132 30)
\fcir f:0 r:1.5
\move(132 33)
\fcir f:0 r:1.5
\move(135 0)
\fcir f:0 r:1.5
\move(135 3)
\fcir f:0 r:1.5
\move(135 6)
\fcir f:0 r:1.5
\move(135 9)
\fcir f:0 r:1.5
\move(135 12)
\fcir f:0 r:1.5
\move(135 15)
\fcir f:0 r:1.5
\move(135 18)
\fcir f:0 r:1.5
\move(135 21)
\fcir f:0 r:1.5
\move(135 24)
\fcir f:0 r:1.5
\move(135 27)
\fcir f:0 r:1.5
\move(135 30)
\fcir f:0 r:1.5
\move(135 33)
\fcir f:0 r:1.5
\move(138 0)
\fcir f:0 r:1.5
\move(138 3)
\fcir f:0 r:1.5
\move(138 6)
\fcir f:0 r:1.5
\move(138 9)
\fcir f:0 r:1.5
\move(138 12)
\fcir f:0 r:1.5
\move(138 15)
\fcir f:0 r:1.5
\move(138 18)
\fcir f:0 r:1.5
\move(138 21)
\fcir f:0 r:1.5
\move(138 24)
\fcir f:0 r:1.5
\move(138 27)
\fcir f:0 r:1.5
\move(138 30)
\fcir f:0 r:1.5
\move(138 33)
\fcir f:0 r:1.5
\move(141 0)
\fcir f:0 r:1.5
\move(141 3)
\fcir f:0 r:1.5
\move(141 6)
\fcir f:0 r:1.5
\move(141 9)
\fcir f:0 r:1.5
\move(141 12)
\fcir f:0 r:1.5
\move(141 15)
\fcir f:0 r:1.5
\move(141 18)
\fcir f:0 r:1.5
\move(141 21)
\fcir f:0 r:1.5
\move(141 24)
\fcir f:0 r:1.5
\move(141 27)
\fcir f:0 r:1.5
\move(141 30)
\fcir f:0 r:1.5
\move(141 33)
\fcir f:0 r:1.5
\move(144 0)
\fcir f:0 r:1.5
\move(144 3)
\fcir f:0 r:1.5
\move(144 6)
\fcir f:0 r:1.5
\move(144 9)
\fcir f:0 r:1.5
\move(144 12)
\fcir f:0 r:1.5
\move(144 15)
\fcir f:0 r:1.5
\move(144 18)
\fcir f:0 r:1.5
\move(144 21)
\fcir f:0 r:1.5
\move(144 24)
\fcir f:0 r:1.5
\move(144 27)
\fcir f:0 r:1.5
\move(144 30)
\fcir f:0 r:1.5
\move(144 33)
\fcir f:0 r:1.5
\move(144 36)
\fcir f:0 r:1.5
\move(147 0)
\fcir f:0 r:1.5
\move(147 3)
\fcir f:0 r:1.5
\move(147 6)
\fcir f:0 r:1.5
\move(147 9)
\fcir f:0 r:1.5
\move(147 12)
\fcir f:0 r:1.5
\move(147 15)
\fcir f:0 r:1.5
\move(147 18)
\fcir f:0 r:1.5
\move(147 21)
\fcir f:0 r:1.5
\move(147 24)
\fcir f:0 r:1.5
\move(147 27)
\fcir f:0 r:1.5
\move(147 30)
\fcir f:0 r:1.5
\move(147 33)
\fcir f:0 r:1.5
\move(147 36)
\fcir f:0 r:1.5
\move(150 0)
\fcir f:0 r:1.5
\move(150 3)
\fcir f:0 r:1.5
\move(150 6)
\fcir f:0 r:1.5
\move(150 9)
\fcir f:0 r:1.5
\move(150 12)
\fcir f:0 r:1.5
\move(150 15)
\fcir f:0 r:1.5
\move(150 18)
\fcir f:0 r:1.5
\move(150 21)
\fcir f:0 r:1.5
\move(150 24)
\fcir f:0 r:1.5
\move(150 27)
\fcir f:0 r:1.5
\move(150 30)
\fcir f:0 r:1.5
\move(150 33)
\fcir f:0 r:1.5
\move(150 36)
\fcir f:0 r:1.5
\move(153 0)
\fcir f:0 r:1.5
\move(153 3)
\fcir f:0 r:1.5
\move(153 6)
\fcir f:0 r:1.5
\move(153 9)
\fcir f:0 r:1.5
\move(153 12)
\fcir f:0 r:1.5
\move(153 15)
\fcir f:0 r:1.5
\move(153 18)
\fcir f:0 r:1.5
\move(153 21)
\fcir f:0 r:1.5
\move(153 24)
\fcir f:0 r:1.5
\move(153 27)
\fcir f:0 r:1.5
\move(153 30)
\fcir f:0 r:1.5
\move(153 33)
\fcir f:0 r:1.5
\move(153 36)
\fcir f:0 r:1.5
\move(156 0)
\fcir f:0 r:1.5
\move(156 3)
\fcir f:0 r:1.5
\move(156 6)
\fcir f:0 r:1.5
\move(156 9)
\fcir f:0 r:1.5
\move(156 12)
\fcir f:0 r:1.5
\move(156 15)
\fcir f:0 r:1.5
\move(156 18)
\fcir f:0 r:1.5
\move(156 21)
\fcir f:0 r:1.5
\move(156 24)
\fcir f:0 r:1.5
\move(156 27)
\fcir f:0 r:1.5
\move(156 30)
\fcir f:0 r:1.5
\move(156 33)
\fcir f:0 r:1.5
\move(156 36)
\fcir f:0 r:1.5
\move(159 0)
\fcir f:0 r:1.5
\move(159 3)
\fcir f:0 r:1.5
\move(159 6)
\fcir f:0 r:1.5
\move(159 9)
\fcir f:0 r:1.5
\move(159 12)
\fcir f:0 r:1.5
\move(159 15)
\fcir f:0 r:1.5
\move(159 18)
\fcir f:0 r:1.5
\move(159 21)
\fcir f:0 r:1.5
\move(159 24)
\fcir f:0 r:1.5
\move(159 27)
\fcir f:0 r:1.5
\move(159 30)
\fcir f:0 r:1.5
\move(159 33)
\fcir f:0 r:1.5
\move(159 36)
\fcir f:0 r:1.5
\move(162 0)
\fcir f:0 r:1.5
\move(162 3)
\fcir f:0 r:1.5
\move(162 6)
\fcir f:0 r:1.5
\move(162 9)
\fcir f:0 r:1.5
\move(162 12)
\fcir f:0 r:1.5
\move(162 15)
\fcir f:0 r:1.5
\move(162 18)
\fcir f:0 r:1.5
\move(162 21)
\fcir f:0 r:1.5
\move(162 24)
\fcir f:0 r:1.5
\move(162 27)
\fcir f:0 r:1.5
\move(162 30)
\fcir f:0 r:1.5
\move(162 33)
\fcir f:0 r:1.5
\move(162 36)
\fcir f:0 r:1.5
\move(165 0)
\fcir f:0 r:1.5
\move(165 3)
\fcir f:0 r:1.5
\move(165 6)
\fcir f:0 r:1.5
\move(165 9)
\fcir f:0 r:1.5
\move(165 12)
\fcir f:0 r:1.5
\move(165 15)
\fcir f:0 r:1.5
\move(165 18)
\fcir f:0 r:1.5
\move(165 21)
\fcir f:0 r:1.5
\move(165 24)
\fcir f:0 r:1.5
\move(165 27)
\fcir f:0 r:1.5
\move(165 30)
\fcir f:0 r:1.5
\move(165 33)
\fcir f:0 r:1.5
\move(165 36)
\fcir f:0 r:1.5
\move(168 0)
\fcir f:0 r:1.5
\move(168 3)
\fcir f:0 r:1.5
\move(168 6)
\fcir f:0 r:1.5
\move(168 9)
\fcir f:0 r:1.5
\move(168 12)
\fcir f:0 r:1.5
\move(168 15)
\fcir f:0 r:1.5
\move(168 18)
\fcir f:0 r:1.5
\move(168 21)
\fcir f:0 r:1.5
\move(168 24)
\fcir f:0 r:1.5
\move(168 27)
\fcir f:0 r:1.5
\move(168 30)
\fcir f:0 r:1.5
\move(168 33)
\fcir f:0 r:1.5
\move(168 36)
\fcir f:0 r:1.5
\move(82 -10)
\lvec(82 40)
\move(94 -10)
\lvec(94 40)
\htext(82 -10){{\scriptsize $D_2$}}
\htext(40 -10){{\scriptsize $D_1$}}
\htext(120 -10){{\scriptsize $D_3$}}
\end{texdraw}
$$
\caption{Division of a diagram}\label{divdiag}
\end{figure}

We are going to show the non-speciality of a large class of systems
not covered by Propositions \ref{easyf1}, \ref{easyfn} and \ref{easyf0}. In fact
we want to construct a finite set $\mathcal E$ of systems and prove
that for $L \notin \mathcal E$ the Conjecture \ref{mainconj} holds.
We will provide conditions on $m$, $n$, $b$ and $a$ such that, under
these conditions, the diagram for $\sys_n(a,b)$, i.e
$\diag([1]^{\times n},\dots,[b]^{\times n},[b+1]^{\times (a+1)})$
can be divided into three parts (see Figure \ref{divdiag})
$D_1+D_2+D_3$ such that $D_2$ is fixed, while $D_1$ and $D_3$ can
be reduced ($D_1$ from the left, $D_3$ from the right) to
much smaller diagrams $D_1'$ and $D_3'$. We will show
that it is possible to enumerate
all reductions of $D_1$ and $D_3$ and show that in each case
the system $\sys(D_1'+D_2+D_3';m^{\times r})$ is non-special.
Then we will apply Theorem \ref{Mred}. 
For example, for $n\geq 4$ and $b\geq 8$ the diagram
can be always written as follows:
$$D_1+\diag([8]^{\times 4})+D_3.$$
It is rather clear that we can $3$-reduce the above to
$$\diag(a_1,a_2)+\diag([8]^{\times 4})+\diag(b_1,b_2),$$
where the set of possible $a_1$, $a_2$, $b_1$, $b_2$ is finite.

\begin{figure}[ht!]
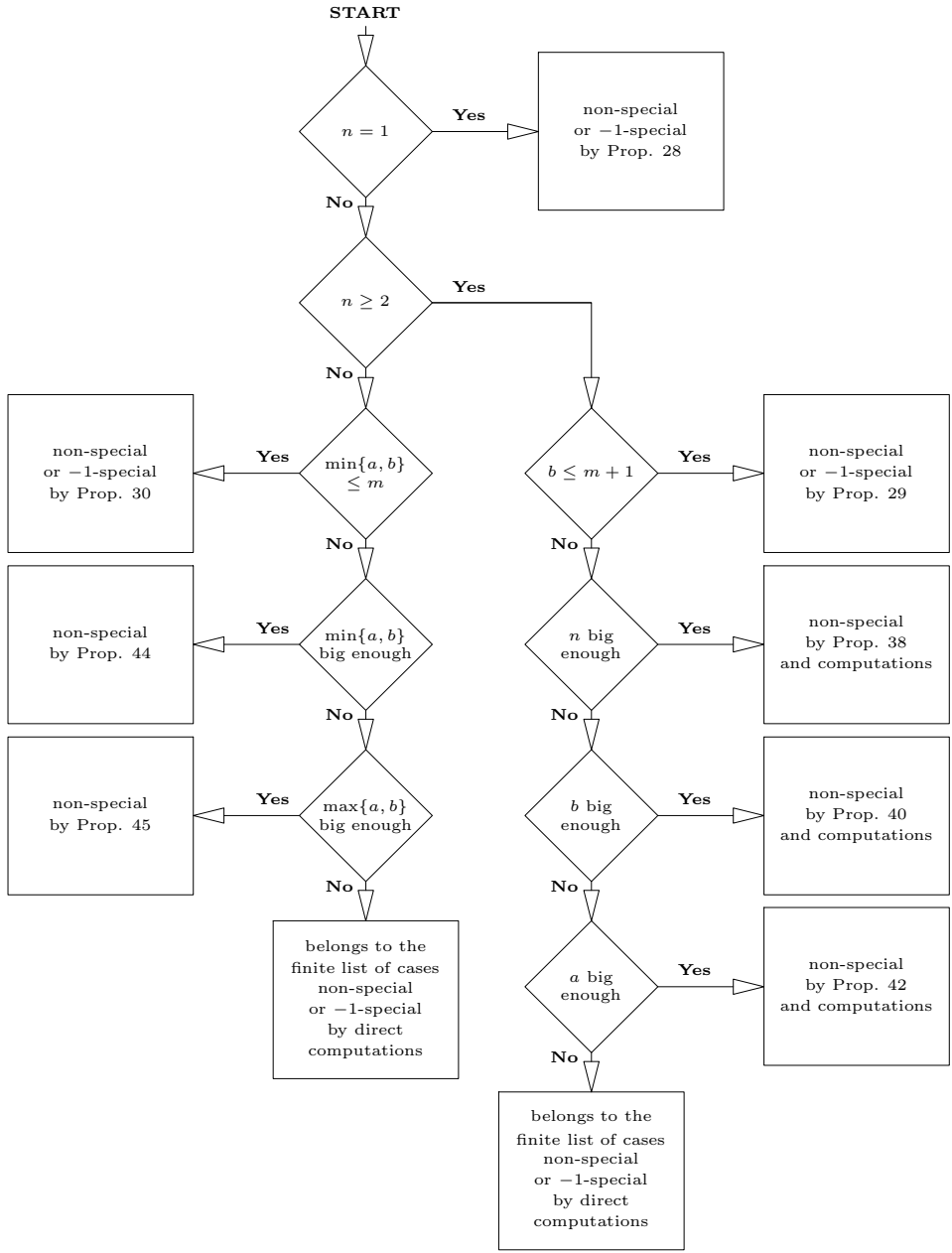

\centertexdraw{
\drawdim pt
\linewd 0.3
\textref h:C v:C
\move(100 15)
\avec(100 0)
\htext(100 20){{\tiny \bf START}}
\move(100 0)
\lvec(125 -25)
\lvec(100 -50)
\lvec(75 -25)
\lvec(100 0)
\htext(100 -25){{\tiny $n=1$}}
\move(100 -50)
\avec(100 -65)
\htext(90 -52){{\tiny \bf No}}
\move(125 -25)
\avec(165 -25)
\htext(139 -19){{\tiny \bf Yes}}
\move(165 5)
\lvec(235 5)
\lvec(235 -55)
\lvec(165 -55)
\lvec(165 5)
\htext(200 -17){{\tiny non-special}}
\htext(200 -25){{\tiny or $-1$-special}}
\htext(200 -33){{\tiny by Prop. \ref{easyf1}}}
\move(100 -65)
\lvec(125 -90)
\lvec(100 -115)
\lvec(75 -90)
\lvec(100 -65)
\htext(100 -90){{\tiny $n \geq 2$}}
\move(100 -115)
\avec(100 -130)
\htext(90 -117){{\tiny \bf No}}
\move(125 -90)
\lvec(185 -90)
\avec(185 -130)
\htext(139 -84){{\tiny \bf Yes}}
\move(100 -130)
\lvec(125 -155)
\lvec(100 -180)
\lvec(75 -155)
\lvec(100 -130)
\htext(100 -152){{\tiny $\min\{a,b\}$}}
\htext(100 -159){{\tiny $\leq m$}}
\move(100 -180)
\avec(100 -195)
\htext(90 -182){{\tiny \bf No}}
\move(75 -155)
\avec(35 -155)
\htext(65 -149){{\tiny \bf Yes}}
\move(185 -130)
\lvec(210 -155)
\lvec(185 -180)
\lvec(160 -155)
\lvec(185 -130)
\htext(185 -155){{\tiny $b \leq m+1$}}
\move(185 -180)
\avec(185 -195)
\htext(175 -182){{\tiny \bf No}}
\move(210 -155)
\avec(250 -155)
\htext(224 -149){{\tiny \bf Yes}}
\move(250 -125)
\lvec(320 -125)
\lvec(320 -185)
\lvec(250 -185)
\lvec(250 -125)
\htext(285 -147){{\tiny non-special}}
\htext(285 -155){{\tiny or $-1$-special}}
\htext(285 -163){{\tiny by Prop. \ref{easyfn}}}
\move(35 -125)
\lvec(-35 -125)
\lvec(-35 -185)
\lvec(35 -185)
\lvec(35 -125)
\htext(0 -147){{\tiny non-special}}
\htext(0 -155){{\tiny or $-1$-special}}
\htext(0 -163){{\tiny by Prop. \ref{easyf0}}}
\move(185 -195)
\lvec(210 -220)
\lvec(185 -245)
\lvec(160 -220)
\lvec(185 -195)
\htext(185 -217){{\tiny $n$ big}}
\htext(185 -224){{\tiny enough}}
\move(185 -245)
\avec(185 -260)
\htext(175 -247){{\tiny \bf No}}
\move(210 -220)
\avec(250 -220)
\htext(224 -214){{\tiny \bf Yes}}
\move(250 -190)
\lvec(320 -190)
\lvec(320 -250)
\lvec(250 -250)
\lvec(250 -190)
\htext(285 -212){{\tiny non-special}}
\htext(285 -220){{\tiny by Prop. \ref{setbignb}}}
\htext(285 -228){{\tiny and computations}}
\move(185 -260)
\lvec(210 -285)
\lvec(185 -310)
\lvec(160 -285)
\lvec(185 -260)
\htext(185 -282){{\tiny $b$ big}}
\htext(185 -289){{\tiny enough}}
\move(185 -310)
\avec(185 -325)
\htext(175 -312){{\tiny \bf No}}
\move(210 -285)
\avec(250 -285)
\htext(224 -279){{\tiny \bf Yes}}
\move(250 -255)
\lvec(320 -255)
\lvec(320 -315)
\lvec(250 -315)
\lvec(250 -255)
\htext(285 -277){{\tiny non-special}}
\htext(285 -285){{\tiny by Prop. \ref{setnb}}}
\htext(285 -293){{\tiny and computations}}
\move(185 -325)
\lvec(210 -350)
\lvec(185 -375)
\lvec(160 -350)
\lvec(185 -325)
\htext(185 -347){{\tiny $a$ big}}
\htext(185 -354){{\tiny enough}}
\move(185 -375)
\avec(185 -390)
\htext(175 -377){{\tiny \bf No}}
\move(210 -350)
\avec(250 -350)
\htext(224 -344){{\tiny \bf Yes}}
\move(250 -320)
\lvec(320 -320)
\lvec(320 -380)
\lvec(250 -380)
\lvec(250 -320)
\htext(285 -342){{\tiny non-special}}
\htext(285 -350){{\tiny by Prop. \ref{setnba}}}
\htext(285 -358){{\tiny and computations}}
\move(150 -390)
\lvec(220 -390)
\lvec(220 -450)
\lvec(150 -450)
\lvec(150 -390)
\htext(185 -400){{\tiny belongs to the}}
\htext(185 -408){{\tiny finite list of cases}}
\htext(185 -416){{\tiny non-special}}
\htext(185 -424){{\tiny or $-1$-special}}
\htext(185 -432){{\tiny by direct}}
\htext(185 -440){{\tiny computations}}
\move(100 -195)
\lvec(125 -220)
\lvec(100 -245)
\lvec(75 -220)
\lvec(100 -195)
\htext(100 -217){{\tiny $\min\{a,b\}$}}
\htext(100 -224){{\tiny big enough}}
\move(100 -245)
\avec(100 -260)
\htext(90 -247){{\tiny \bf No}}
\move(75 -220)
\avec(35 -220)
\htext(65 -214){{\tiny \bf Yes}}
\move(35 -190)
\lvec(-35 -190)
\lvec(-35 -250)
\lvec(35 -250)
\lvec(35 -190)
\htext(0 -216){{\tiny non-special}}
\htext(0 -224){{\tiny by Prop. \ref{setpb}}}
\move(100 -260)
\lvec(125 -285)
\lvec(100 -310)
\lvec(75 -285)
\lvec(100 -260)
\htext(100 -282){{\tiny $\max\{a,b\}$}}
\htext(100 -289){{\tiny big enough}}
\move(100 -310)
\avec(100 -325)
\htext(90 -312){{\tiny \bf No}}
\move(75 -285)
\avec(35 -285)
\htext(65 -279){{\tiny \bf Yes}}
\move(35 -255)
\lvec(-35 -255)
\lvec(-35 -315)
\lvec(35 -315)
\lvec(35 -255)
\htext(0 -281){{\tiny non-special}}
\htext(0 -289){{\tiny by Prop. \ref{setpba}}}
\move(65 -325)
\lvec(135 -325)
\lvec(135 -385)
\lvec(65 -385)
\lvec(65 -325)
\htext(100 -335){{\tiny belongs to the}}
\htext(100 -343){{\tiny finite list of cases}}
\htext(100 -351){{\tiny non-special}}
\htext(100 -359){{\tiny or $-1$-special}}
\htext(100 -367){{\tiny by direct}}
\htext(100 -375){{\tiny computations}}
}
\caption{The way of proving the main conjecture}\label{idea}
\end{figure}

To make the understanding of our proof easier, we present the outline
on Figure \ref{idea}.

Define (for a diagram $D$ and $m \geq 2$) the following number:
$$p(D) = \left\lfloor \frac{\# D}{\binom{m+1}{2}} \right\rfloor.$$
Observe that if $\sys(D;m^{\times p(D)})$ and $\sys(D;m^{\times (p(D)+1)})$
are non-special then $\sys(D;m^{\times r})$ is non-special for $r \geq 0$.

For $D=\diag(a_1,\dots,a_s)$ let 
$\rev(D) = \diag(a_s,\dots,a_1)$.

\begin{proposition}
\label{setbignb}
Let $m$, $N$, $B$ be integers, $m \geq 2$, $N \geq m$, $B \geq m+2$.
There exists the finite set $\mathcal D$ of diagrams, such that
if for all $D \in \mathcal D$ both systems $\sys(D;m^{\times p(D)})$
and $\sys(D;m^{\times (p(D)+1)})$ are non-special then
for any $n \geq N$, $b \geq B$, $a \geq 0$, $r \geq 0$
the system $\sys_n(a,b;m^{\times r})$ is non-special.
Moreover, the set $\mathcal D$ can be found algorithmically.
\end{proposition}

\begin{proof}
The idea is to choose $\mathcal D$ to be the set of diagrams with the following
property:
any diagram for $\sys_n(a,b;m^{\times r})$, i.e. diagram
$$D=\diag([1]^{\times n},[2]^{\times n},\dots,[b]^{\times n},b+1,[b+1]^{\times a})$$
can be reduced $r$ times, or can be reduced to a diagram from $\mathcal D$.
If this is the case then we conclude by Theorem \ref{Mred}.
Of course, we use $m$-reductions.
We will show how $D$ can be reduced, and simultanously we will construct
$\mathcal D$.

We begin with $m$-reducing from the left, the first layer being the lowest one
(see Figure \ref{reds1}).
\begin{figure}[ht!]
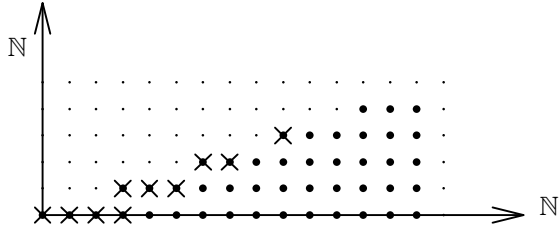

$$
\begin{texdraw}
\drawdim pt
\move(0 0)
\fcir f:0 r:0.5
\move(0 10)
\fcir f:0 r:0.5
\move(0 20)
\fcir f:0 r:0.5
\move(0 30)
\fcir f:0 r:0.5
\move(0 40)
\fcir f:0 r:0.5
\move(0 50)
\fcir f:0 r:0.5
\move(10 0)
\fcir f:0 r:0.5
\move(10 10)
\fcir f:0 r:0.5
\move(10 20)
\fcir f:0 r:0.5
\move(10 30)
\fcir f:0 r:0.5
\move(10 40)
\fcir f:0 r:0.5
\move(10 50)
\fcir f:0 r:0.5
\move(20 0)
\fcir f:0 r:0.5
\move(20 10)
\fcir f:0 r:0.5
\move(20 20)
\fcir f:0 r:0.5
\move(20 30)
\fcir f:0 r:0.5
\move(20 40)
\fcir f:0 r:0.5
\move(20 50)
\fcir f:0 r:0.5
\move(30 0)
\fcir f:0 r:0.5
\move(30 10)
\fcir f:0 r:0.5
\move(30 20)
\fcir f:0 r:0.5
\move(30 30)
\fcir f:0 r:0.5
\move(30 40)
\fcir f:0 r:0.5
\move(30 50)
\fcir f:0 r:0.5
\move(40 0)
\fcir f:0 r:0.5
\move(40 10)
\fcir f:0 r:0.5
\move(40 20)
\fcir f:0 r:0.5
\move(40 30)
\fcir f:0 r:0.5
\move(40 40)
\fcir f:0 r:0.5
\move(40 50)
\fcir f:0 r:0.5
\move(50 0)
\fcir f:0 r:0.5
\move(50 10)
\fcir f:0 r:0.5
\move(50 20)
\fcir f:0 r:0.5
\move(50 30)
\fcir f:0 r:0.5
\move(50 40)
\fcir f:0 r:0.5
\move(50 50)
\fcir f:0 r:0.5
\move(60 0)
\fcir f:0 r:0.5
\move(60 10)
\fcir f:0 r:0.5
\move(60 20)
\fcir f:0 r:0.5
\move(60 30)
\fcir f:0 r:0.5
\move(60 40)
\fcir f:0 r:0.5
\move(60 50)
\fcir f:0 r:0.5
\move(70 0)
\fcir f:0 r:0.5
\move(70 10)
\fcir f:0 r:0.5
\move(70 20)
\fcir f:0 r:0.5
\move(70 30)
\fcir f:0 r:0.5
\move(70 40)
\fcir f:0 r:0.5
\move(70 50)
\fcir f:0 r:0.5
\move(80 0)
\fcir f:0 r:0.5
\move(80 10)
\fcir f:0 r:0.5
\move(80 20)
\fcir f:0 r:0.5
\move(80 30)
\fcir f:0 r:0.5
\move(80 40)
\fcir f:0 r:0.5
\move(80 50)
\fcir f:0 r:0.5
\move(90 0)
\fcir f:0 r:0.5
\move(90 10)
\fcir f:0 r:0.5
\move(90 20)
\fcir f:0 r:0.5
\move(90 30)
\fcir f:0 r:0.5
\move(90 40)
\fcir f:0 r:0.5
\move(90 50)
\fcir f:0 r:0.5
\move(100 0)
\fcir f:0 r:0.5
\move(100 10)
\fcir f:0 r:0.5
\move(100 20)
\fcir f:0 r:0.5
\move(100 30)
\fcir f:0 r:0.5
\move(100 40)
\fcir f:0 r:0.5
\move(100 50)
\fcir f:0 r:0.5
\move(110 0)
\fcir f:0 r:0.5
\move(110 10)
\fcir f:0 r:0.5
\move(110 20)
\fcir f:0 r:0.5
\move(110 30)
\fcir f:0 r:0.5
\move(110 40)
\fcir f:0 r:0.5
\move(110 50)
\fcir f:0 r:0.5
\move(120 0)
\fcir f:0 r:0.5
\move(120 10)
\fcir f:0 r:0.5
\move(120 20)
\fcir f:0 r:0.5
\move(120 30)
\fcir f:0 r:0.5
\move(120 40)
\fcir f:0 r:0.5
\move(120 50)
\fcir f:0 r:0.5
\move(130 0)
\fcir f:0 r:0.5
\move(130 10)
\fcir f:0 r:0.5
\move(130 20)
\fcir f:0 r:0.5
\move(130 30)
\fcir f:0 r:0.5
\move(130 40)
\fcir f:0 r:0.5
\move(130 50)
\fcir f:0 r:0.5
\move(140 0)
\fcir f:0 r:0.5
\move(140 10)
\fcir f:0 r:0.5
\move(140 20)
\fcir f:0 r:0.5
\move(140 30)
\fcir f:0 r:0.5
\move(140 40)
\fcir f:0 r:0.5
\move(140 50)
\fcir f:0 r:0.5
\move(150 0)
\fcir f:0 r:0.5
\move(150 10)
\fcir f:0 r:0.5
\move(150 20)
\fcir f:0 r:0.5
\move(150 30)
\fcir f:0 r:0.5
\move(150 40)
\fcir f:0 r:0.5
\move(150 50)
\fcir f:0 r:0.5
\arrowheadtype t:V
\move(0 0)
\avec(180 0)
\move(0 0)
\avec(0 80)
\htext(186 0){$\mathbb{N}$}
\htext(-13 60){$\mathbb{N}$}
\move(0 0)
\fcir f:0 r:1.5
\move(10 0)
\fcir f:0 r:1.5
\move(20 0)
\fcir f:0 r:1.5
\move(30 0)
\fcir f:0 r:1.5
\move(30 10)
\fcir f:0 r:1.5
\move(40 0)
\fcir f:0 r:1.5
\move(40 10)
\fcir f:0 r:1.5
\move(50 0)
\fcir f:0 r:1.5
\move(50 10)
\fcir f:0 r:1.5
\move(60 0)
\fcir f:0 r:1.5
\move(60 10)
\fcir f:0 r:1.5
\move(60 20)
\fcir f:0 r:1.5
\move(70 0)
\fcir f:0 r:1.5
\move(70 10)
\fcir f:0 r:1.5
\move(70 20)
\fcir f:0 r:1.5
\move(80 0)
\fcir f:0 r:1.5
\move(80 10)
\fcir f:0 r:1.5
\move(80 20)
\fcir f:0 r:1.5
\move(90 0)
\fcir f:0 r:1.5
\move(90 10)
\fcir f:0 r:1.5
\move(90 20)
\fcir f:0 r:1.5
\move(90 30)
\fcir f:0 r:1.5
\move(100 0)
\fcir f:0 r:1.5
\move(100 10)
\fcir f:0 r:1.5
\move(100 20)
\fcir f:0 r:1.5
\move(100 30)
\fcir f:0 r:1.5
\move(110 0)
\fcir f:0 r:1.5
\move(110 10)
\fcir f:0 r:1.5
\move(110 20)
\fcir f:0 r:1.5
\move(110 30)
\fcir f:0 r:1.5
\move(120 0)
\fcir f:0 r:1.5
\move(120 10)
\fcir f:0 r:1.5
\move(120 20)
\fcir f:0 r:1.5
\move(120 30)
\fcir f:0 r:1.5
\move(120 40)
\fcir f:0 r:1.5
\move(130 0)
\fcir f:0 r:1.5
\move(130 10)
\fcir f:0 r:1.5
\move(130 20)
\fcir f:0 r:1.5
\move(130 30)
\fcir f:0 r:1.5
\move(130 40)
\fcir f:0 r:1.5
\move(140 0)
\fcir f:0 r:1.5
\move(140 10)
\fcir f:0 r:1.5
\move(140 20)
\fcir f:0 r:1.5
\move(140 30)
\fcir f:0 r:1.5
\move(140 40)
\fcir f:0 r:1.5
\move(0 0)
\move(-3 -3)
\lvec(3 3)
\move(3 -3)
\lvec(-3 3)
\move(10 0)
\move(7 -3)
\lvec(13 3)
\move(13 -3)
\lvec(7 3)
\move(20 0)
\move(17 -3)
\lvec(23 3)
\move(23 -3)
\lvec(17 3)
\move(30 10)
\move(27 7)
\lvec(33 13)
\move(33 7)
\lvec(27 13)
\move(30 0)
\move(27 -3)
\lvec(33 3)
\move(33 -3)
\lvec(27 3)
\move(40 10)
\move(37 7)
\lvec(43 13)
\move(43 7)
\lvec(37 13)
\move(50 10)
\move(47 7)
\lvec(53 13)
\move(53 7)
\lvec(47 13)
\move(60 20)
\move(57 17)
\lvec(63 23)
\move(63 17)
\lvec(57 23)
\move(70 20)
\move(67 17)
\lvec(73 23)
\move(73 17)
\lvec(67 23)
\move(90 30)
\move(87 27)
\lvec(93 33)
\move(93 27)
\lvec(87 33)
\end{texdraw}
$$
\caption{Reduction from the left}\label{reds1}
\end{figure}
After performing $n$ such reductions, we obtain
$$D_6=\diag([m+1]^{\times n},[m+2]^{\times n},\dots,[b]^{\times n},[b+1]^{\times (a+1)}).$$

We must deal with two cases, $B \geq 2m-1$ and $B < 2m-1$, separately.

{\bf Case $B \geq 2m-1$.} For $b \geq B \geq 2m-1$ we can write
$$D_6=\diag([m+1]^{\times n},\dots,[2m-2]^{\times n},[2m-1]^{\times n},\dots,[b+1]^{\times (a+1)}).$$
We will show that this diagram can be reduced to 
$$D_5=\diag([m+1]^{\times n},\dots,[2m-2]^{\times n},a_1,\dots,a_{m-1}),$$
where $(a_1,\dots,a_{m-1})$ is an admissible $([2m-1]^{\times m})$-tail.
To see the above, take
\begin{align*}
G & =\diag([2m-1]^{\times n},\dots,[b+1]^{\times (a+1)}) \\
  & =\diag([2m-1]^{\times m},[2m-1]^{\times (n-m)},\dots,[b+1]^{\times (a+1)})
\end{align*}
and use Proposition \ref{condonred}.
It means, in particular, that $\diag(a_1,\dots,a_{m-1})$ can be obtained from
$\diag([2m-1]^{\times m},[\symbx]^{\times (m-1)})$ by a sequence of
symbolic reductions. Put 
$$\mathcal D_5=\{ G : G \text{ is an admissible $([2m-1]^{\times m})$-tail}\}.$$

Now take
$$
D_4 = \rev(D_5) = \diag(a_{m-1},\dots,a_1,[2m-2]^{\times n},\dots,[m+1]^{\times n}).
$$
$D_4$ can be reduced to a diagram
$$D_3 = \diag(a_{m-1},\dots,a_1,[2m-2]^{\times n},\dots,[m+2]^{\times n},b_1,\dots,b_{m-1}),$$
for an admissible $(m+1)$-$([0]^{\times (m-1)})$-tail $(b_1,\dots,b_{m-1})$.
Taking
$$\mathcal D_3=\{\text{all admissible $(m+1)$-$([0]^{\times (m-1)})$-tails}\}$$
we will have
$$D_3 \in \{\rev(G)+\diag([2m-2]^{\times n},\dots,[m+2]^{\times n})+H : G \in \mathcal D_5, \, H \in \mathcal D_3\}.$$
Again, $D_3$ can be reduced to 
$$D_2 = \diag(a_{m-1},\dots,a_1,[2m-2]^{\times n},\dots,[m+3]^{\times n},c_1,\dots,c_{m-1}),$$
for an admissible $(m+2)$-$(b_1,\dots,b_{m-1})$-tail $(c_1,\dots,c_{m-1})$.
Taking
$$\mathcal D_2=\{\text{all admissible $(m+2)$-$G$-tails } : G \in \mathcal D_3\}$$
we will have
$$D_2 \in \{\rev(G)+\diag([2m-2]^{\times n},\dots,[m+3]^{\times n})+H : G \in \mathcal D_5, \, H \in \mathcal D_2\}.$$
This can be repeated for each $j=m+3,\dots,2m-3$ until the following is obtained:
\begin{align*}
D_1 & = \diag(a_{m-1},\dots,a_1,[2m-2]^{\times n},d_1,\dots,d_{m-1}),\\
\mathcal D_1 & = \{\text{all admissible $(2m-3)$-$G$-tails } : G \text{ in the previous set } \mathcal D\},\\
D_1 & \in \{\rev(G)+\diag([2m-2]^{\times n})+H : G \in \mathcal D_5, \, H \in \mathcal D_1\}.
\end{align*}
Now we do the above once more for $2m-2$, but leaving the
part $\diag([2m-2]^{\times N})$ untouched, in order to finish with a diagram
big enough to obtain non-speciality.
So we obtain
\begin{align*}
D_0 & = \diag(a_{m-1},\dots,a_1,[2m-2]^{\times N},e_1,\dots,e_{m-1}),\\
\mathcal D_0 & = \{\text{all admissible $(2m-2)$-$G$-tails } : G \in \mathcal D_1\},\\
D_0 & \in \{\rev(G)+\diag([2m-2]^{\times N})+H : G \in \mathcal D_5, \, H \in \mathcal D_0\}.
\end{align*}
Since $D_0$ has been obtained from $D$ by a sequence of $m$-reductions, putting
$$
\mathcal D = \{\rev(G)+\diag([2m-2]^{\times N})+H : G \in \mathcal D_5, \, H \in \mathcal D_0\}.
$$
we are done.

{\bf Case $B < 2m-1$.} For each $b \in \{B, \dots, 2m-2\}$ we do the following.
Put
$$D_6=\diag([m+1]^{\times n},\dots,[b+1]^{\times (a+1)}).$$
By Proposition \ref{condonred}, this diagram can be reduced to 
$$D_5=\diag([m+1]^{\times n},\dots,[b]^{\times n},a_1,\dots,a_{m-1}),$$
where $(a_1,\dots,a_{m-1})$ is an admissible $(b+1)$-$([0]^{\times (m-1)})$-tail.
Put 
$$\mathcal D_5=\{ G : G \text{ is an admissible $(b+1)$-$([0]^{\times (m-1)})$-tail}\}.$$

Now take
$$
D_4 = \rev(D_5) = \diag(a_{m-1},\dots,a_1,[b]^{\times n},\dots,[m+1]^{\times n}).
$$
As in the previous case, we repeat reducing together with generating all
admissible $j$-$G$-tails, until the following is obtained:
\begin{align*}
D_1 & = \diag(a_{m-1},\dots,a_1,[b]^{\times n},d_1,\dots,d_{m-1}),\\
\mathcal D_1 & = \{\text{all admissible $(b-1)$-$G$-tails } : G \text{ in the previous set } \mathcal D\},\\
D_1 & \in \{\rev(G)+\diag([b]^{\times n})+H : G \in \mathcal D_5, \, H \in \mathcal D_1\}.
\end{align*}
Now, as before, we do the above once more for $b$, but leaving the
part $\diag([b]^{\times N})$ untouched.
\end{proof}

\begin{example}
Let us show how we reduce for $(m,N,B) = (4,10,7)$
and $(m,N,B) = (4,12,6)$. In the first case
we consider
$m=4$, $n\geq 10$, $b \geq 7$, $a \geq 0$ and take
$$D=\diag([1]^{\times n},\dots,[b]^{\times n},[b+1]^{\times (a+1)}).$$
In order to make reductions, we consider
\begin{multline*}
D = \diag([1]^{\times n},\dots,[4]^{\times n}) + \diag([5]^{\times n})
+\diag([6]^{\times (n-10)}) + \diag([6]^{\times 10}) \\
+ \diag(7,7,7,7,\symbx,\symbx,\symbx).
\end{multline*}
The diagram will be reduced from left and right, without touching
$\diag([6]^{\times 10})$.

In the second case
we take
$m=4$, $n\geq 12$, $b = 6$, $a \geq 0$,
$$D=\diag([1]^{\times n},\dots,[6]^{\times n},[7]^{\times (a+1)}).$$
In order to make reductions, we consider
\begin{multline*}
D = \diag([1]^{\times n},\dots,[4]^{\times n}) + \diag([5]^{\times n})
+\diag([6]^{\times (n-12)}) + \diag([6]^{\times 12}) \\
+ \diag([7]^{\times (a+1)}).
\end{multline*}
Again, $\diag([6]^{\times 12})$ remains untouched during reductions.
\end{example}

\begin{proposition}
\label{setnb}
Let $m$, $n$, $B$ be integers, $m \geq 2$, $n \geq 2$, $B \geq 2m-1$.
There exists the finite set $\mathcal D$ of diagrams, such that
if for all $D \in \mathcal D$ both systems $\sys(D;m^{\times p(D)})$
and $\sys(D;m^{\times (p(D)+1)})$ are non-special then
for any $b \geq B$, $a \geq 0$, $r \geq 0$
the system $\sys_n(a,b;m^{\times r})$ is non-special.
Moreover, the set $\mathcal D$ can be found algorithmically.
\end{proposition}

\begin{proof}
Take $G=\diag([1]^{\times n},\dots,[B]^{\times n},B+1)$, let
$k=nB+1-m$.
Write
$$G = \diag(e_1,\dots,e_k) + \diag(d_1,\dots,d_m).$$
Let
$$\mathcal D = \{ \diag(e_1,\dots,e_k) + E : E \text{ is an admissible $(d_1,\dots,d_m)$-tail}\}.$$
Now take $b$, $a$ and $r$ as above, let
$$D = \diag([1]^{\times n},\dots,[b]^{\times n},[b+1]^{\times (a+1)}) = G + \diag(\dots).$$
Since $b\geq 2m-1$ we can use Proposition \ref{condonred} to show that
$D$ can be reduced to the diagram
$$D_1 = \diag(e_1,\dots,e_k) + \diag(a_1,\dots,a_{m-1}),$$
where $(a_1,\dots,a_{m-1})$ is an admissible $(d_1,\dots,d_m)$-tail.
Observe that $D_1 \in \mathcal D$
and conclude with Theorem \ref{Mred}.
\end{proof}

\begin{example}
Let $m=4$, $n=4$, $B=7$. Take $b \geq 7$, $a\geq 0$ and consider
$$D=\diag([1]^{\times 4},\dots,[b]^{\times 4},[b+1]^{\times (a+1)}),$$
which can be written as
$$\diag([1]^{\times 4},\dots,[6]^{\times 4},7)+\diag(7,7,7,8,\symbx,\symbx,\symbx) + \diag(\dots).$$
Now the left hand side remains untouched, while the right hand side will be
reduced.
\end{example}

\begin{proposition}
\label{setnba}
Let $m$, $n$, $b$, $A$ be integers, $m \geq 2$, $n \geq 2$, $b \geq m$, $A \geq 0$.
There exists the finite set $\mathcal D$ of diagrams, such that
if for all $D \in \mathcal D$ both systems $\sys(D;m^{\times p(D)})$
and $\sys(D;m^{\times (p(D)+1)})$ are non-special then
for any $a \geq A$, $r \geq 0$
the system $\sys_n(a,b;m^{\times r})$ is non-special.
Moreover, the set $\mathcal D$ can be found algorithmically.
\end{proposition}

\begin{proof}
Let
\begin{multline*}
\mathcal D = \{ \diag([1]^{\times n},\dots,[b]^{\times n},[b+1]^{\times (A+1)}) + G \\
: G \text{ is an admissible $(b+1)$-$([0]^{\times (m-1)})$-tail}\}.
\end{multline*}
Now take $a$ and $r$ as above, let
$$D = \diag([1]^{\times n},\dots,[b]^{\times n},[b+1]^{\times (a+1)}).$$
Since $b\geq m$ we can use Proposition \ref{condonred} to show that
$D$ can be reduced to the diagram
$$D_1 = \diag([1]^{\times n},\dots,[b]^{\times n},[b+1]^{\times (A+1)},d_1,\dots,d_{m-1}),$$
where $(d_1,\dots,d_{m-1})$ is an admissible $(b+1)$-$([0]^{\times (m-1)})$-tail.
Observe that $D_1 \in \mathcal D$
and conclude using Theorem \ref{Mred}.
\end{proof}

\begin{example}
Let $m=4$, $n=3$, $b=6$, $A=1$. Take $a\geq 1$ and consider
$$D=\diag([1]^{\times 3},\dots,[6]^{\times 3},[7]^{\times (a+1)}),$$
which can be written as
$$\diag([1]^{\times 3},\dots,[6]^{\times 3},7,7)+\diag([7]^{\times (a-1)}).$$
Now the left hand side remains untouched, while the right hand side will be
reduced.
\end{example}

\begin{proposition}
\label{setpb}
Let $m$, $B$ be integers, $m \geq 2$, $B \geq 3(m-1)$.
There exists the finite set $\mathcal D$ of diagrams, such that
if for all $D \in \mathcal D$ both systems $\sys(D;m^{\times p(D)})$
and $\sys(D;m^{\times (p(D)+1)})$ are non-special then
for any $a \geq b \geq B$, $r \geq 0$
the system $\sys_0(a,b;m^{\times r})$ is non-special.
Moreover, the set $\mathcal D$ can be found algorithmically.
\end{proposition}

\begin{proof}
It is enough to take
\begin{multline*}
\mathcal D = \{ \diag(1,\dots,B-m+1) + G \\
: G \text{ is an admissible $(B-m+2,\dots,B+1)$-tail}\}.
\end{multline*}
Indeed, taking $b$, $a$ and $r$ as above, let
\begin{multline*}
D = \diag(1,\dots,b,[b+1]^{\times (a-b+1)},b^{\up 1},\dots,1^{\up b}) \\
= \diag(1,\dots,B-m+1) + \diag(B-m+2,\dots,B+1,\dots)
\end{multline*}
be the diagram for $\sys_0(a,b)$ (see Proposition \ref{diagforp1p1}).
Observe that reducing $(b^{\up 1},\dots,1^{\up b})$ is equivalent to reducing
$(b,\dots,1)$, and this part of the diagram will surely be reduced.
Again, by Proposition \ref{condonred},
$D$ can be reduced to a diagram
belonging to $\mathcal D$, since, by assumption,
$B-m+2 \geq 2m-1$.
\end{proof}

\begin{proposition}
\label{setpba}
Let $m$, $b$, $A$ be integers, $A \geq b \geq m \geq 2$.
There exists the finite set $\mathcal D$ of diagrams, such that
if for all $D \in \mathcal D$ both systems $\sys(D;m^{\times p(D)})$
and $\sys(D;m^{\times (p(D)+1)})$ are non-special then
for any $a \geq A$, $r \geq 0$
the system $\sys_0(a,b;m^{\times r})$ is non-special.
Moreover, the set $\mathcal D$ can be found algorithmically.
\end{proposition}

\begin{proof}
It is enough to take
$$\mathcal D = \{ \diag([b+1]^{\times (A+1)}) + G : G \text{ is an admissible $(b+1)$-$([0]^{\times (m-1)})$-tail}\}.$$
Indeed, take $a$ and $r$ as above and let
$$D = \diag([b+1]^{\times (a+1)}) = \diag([b+1]^{\times (A+1)}) + \diag(b+1,b+1,\dots).$$
Again, by Proposition \ref{condonred},
$D$ can be reduced to a diagram
belonging to $\mathcal D$.
\end{proof}

\begin{proposition}
\label{computed1}
For the following values of $m$, $N$, and $B$
the set $\mathcal D$ from Proposition \ref{setbignb} contains only non-special diagrams:
$$
\begin{array}{c||c|c|c|c|c|c|c}
m                      &  2 &  3 &  4 &  5 &  6 &  7 &  8 \\
\hline
N \text{ for } B > m+2 &  2 &  5 & 11 & 11 & 22 & 25 & 41 \\
N \text{ for } B = m+2 &  2 &  5 & 16 & 30 & 51 & 85 &127
\end{array}
$$
\end{proposition}

\begin{proof}
The proof was completed using suitable computer programs. First, one has to
create the set $\mathcal D$. Next, for each diagram $D \in \mathcal D$
a computation of the rank of two interpolation matrices (for $p(D)$ and $p(D)+1$
points of multiplicity $m$) shows that $D$ is non-special.
All programs can be downloaded from \cite{MYWWW}, together with
files containing the results of running them by the author.
\end{proof}

\begin{proposition}
\label{computed2}
For the following values of $m$, $n$, and $B$
the set $\mathcal D$ from Proposition \ref{setnb}
contains only non-special diagrams:
\begin{align*}
&
\begin{array}{c||c||c|c|c||c|c|c}
m  &     3 & 4 & 4 &          4 & 5  & 5  & 5         \\
\hline
n  & 2,3,4 & 2 & 3 & 4,\dots,10 & 2  & 3  & 4,\dots,10 \\
\hline
B  &     6 & 9 & 8 &          9 & 11 & 10 & 9         
\end{array}
\\
&
\begin{array}{c||c|c|c||c|c||c|c}
m  & 6  & 6  & 6          & 7  & 7          & 8  & 8          \\
\hline
n  & 2  & 3  & 4,\dots,21 & 2  & 3,\dots,24 & 2  & 3,\dots,40 \\
\hline
B  & 13 & 12 & 11         & 15 & 13         & 17 & 15
\end{array}
\end{align*}
\end{proposition}

\begin{proof}
Again we use suitable computer programs.
\end{proof}

Observe that only finite number of triples $(m,n,b)$ satisfying
$$2 \leq m \leq 8, \quad n \geq 2, \quad b \geq m+2$$
are not covered by the two previous Propositions. We will not list all of them,
but only present the number of them:
$$
\begin{array}{r||c|c|c|c|c|c|c}
m & 2 & 3 & 4 & 5 & 6 & 7 & 8 \\
\hline
\# \text{ of triples} & 0  & 3  & 17  &  40 & 92 & 154  & 321
\end{array}
$$
Again using suitable
computer programs we proved the following:

\begin{proposition}
\label{computed3}
For every triple $(m,n,b)$ not covered by Proposition \ref{computed1}
or Proposition \ref{computed2} there exists $A$ such that
the set $\mathcal D$ from Proposition \ref{setnba} contains only non-special
diagrams.
The greatest value of $A$ is shown in the table below:
$$
\begin{array}{r||c|c|c|c|c|c|c}
m & 2 & 3 & 4 & 5 & 6 & 7 & 8 \\
\hline
\max A & 0  & 0  & 1 & 3  & 13  & 22  & 33 
\end{array}
$$
\end{proposition}

\begin{proposition}
\label{computed4}
For the following values of $m$ and $B$
the set $\mathcal D$ from Proposition \ref{setpb} contains only non-special diagrams:
$$
\begin{array}{c||c|c|c|c|c|c|c}
m                      &  2 &  3 &  4 &  5 &  6 &  7 &  8 \\
\hline
B  &  6 &  9 & 12 & 15 & 19 & 21 & 24 \\
\end{array}
$$
\end{proposition}

\begin{proof}
Again we use suitable computer programs. Observe that $B=3m$ is sufficient
for all checked cases except for $m=6$. There are no geometrical explantion
to this fact (all systems $\sys_0(a,b;6^{\times r})$ with $a,b \geq 18$
are non-special), but the system
$\sys(D;6^{\times 6})$ for
$$D = \diag(1,2,3,4,5,6,7,8,9,10,11,12,13,12,11,10,2)$$
is special.
\end{proof}

\begin{proposition}
\label{computed5}
For the following values of $m$, $b$, and $A$
the set $\mathcal D$ from Proposition \ref{setpba}
contains only non-special diagrams:
\begin{align*}
&
\begin{array}{c||c||c|c||c|c||c|c|c}
m  &     2 & 3 & 3         & 4  & 4          & 5  & 5  & 5          \\
\hline
b  & 3,4,5 & 4 & 5,\dots,8 & 5  & 6,\dots,11 & 6  & 7  & 8,\dots,14 \\
\hline
A  &     b & 6 & b         & 15 & b          & 28 & 10 & b         
\end{array}
\\
&
\begin{array}{c||c|c|c|c||c||c|c|c|c}
m  & 6  & 6  & 6  & 6           & 7  & 7  & 7  & 7  & 7         \\
\hline
b  & 7  & 8  & 9  & 10,\dots,18 & 8  & 9  & 10 & 11 & 12,\dots,20 \\
\hline
A  & 50 & 21 & 11 & b           & 84 & 31 & 18 & 12 & b
\end{array}
\\
&
\begin{array}{c||c|c|c|c|c}
m  & 8   & 8  & 8  & 8  & 8           \\
\hline
b  & 9   & 10 & 11 & 12 & 13,\dots,23 \\
\hline
A  & 126 & 43 & 27 & 15 & b
\end{array}
\end{align*}
\end{proposition}

\begin{proof}
Again we use suitable computer programs.
\end{proof}

\section{Final cases}

There are some quadruples $(m,n,a,b)$ not covered by previous Propositions.
For each of these we must find
$r_1$ and $r_2$ such that for all $r \leq r_1$ and all $r \geq r_2$ the system
$\sys_n(a,b;m^{\times r})$ is non-special. This can be done
by direct computations. If $r_1+1<r_2$ then we must check if all systems
for $r_1 < r < r_2$ are $-1$-special.
This was done 
by a computer program. Here we present the number of final $(m,n,a,b)$'s
together with the number of special systems found,
depending on $m$:
$$
\begin{array}{r||c|c|c|c|c|c|c}
m & 2 & 3 & 4 & 5 & 6 & 7 & 8 \\
\hline
\text{number of cases} & 0 & 2 & 11 & 30 & 90 & 187 & 353 \\
\hline
\text{number of special systems} & 0 & 1 & 5 & 12 & 37 & 70 & 134
\end{array}
$$
While checking $-1$-speciality we considered only $-1$-systems
with imposed base points in general position, i.e.
$\sys(d;m_1,\dots,m_r,\overline{k_1,\dots,k_s})$ satisfying
$$\dim \sys(d;m_1,\dots,m_r,\overline{k_1,\dots,k_s})
= \dim \sys(d;m_1,\dots,m_r,k_1,\dots,k_s).$$
Also in Propositions \ref{easyf1}, \ref{easyfn} and \ref{easyf0} we used
$-1$-systems with imposed base points in general position.

\end{document}